\documentclass[11pt]{article}
\usepackage{color}

\usepackage[utf8]{inputenc}
\usepackage[T1]{fontenc}
\usepackage{amsmath,amssymb,epsfig,bbm}
\usepackage{comment}
\usepackage[colorlinks]{hyperref}
\usepackage{todonotes}

\newtheorem{defi}{Definition}
\newtheorem{prop}[defi]{Proposition}
\newtheorem{theo}[defi]{Theorem}
\newtheorem{conj}[defi]{Conjecture}
\newtheorem{lemm}[defi]{Lemma}
\newtheorem{coro}[defi]{Corollary}
\newtheorem{rema}[defi]{Remark}
\newtheorem{exem}[defi]{Example}
\newtheorem{exems}[defi]{Examples}

\newcommand{\bdefi}{\begin{defi}}
\newcommand{\edefi}{\end{defi}}
\newcommand{\bprop}{\begin{prop}}
\newcommand{\eprop}{\end{prop}}
\newcommand{\btheo}{\begin{theo}}
\newcommand{\etheo}{\end{theo}}
\newcommand{\blemm}{\begin{lemm}}
\newcommand{\brema}{\begin{rema}}
\newcommand{\erema}{\end{rema}}
\newcommand{\bexer}{\begin{exem}}
\newcommand{\eexer}{\end{exem}}
\newcommand{\bexems}{\begin{exems}}
\newcommand{\eexems}{\end{exems}}
\newcommand{\bconj}{\begin{conj}}
\newcommand{\econj}{\end{conj}}
\newcommand{\elemm}{\end{lemm}}
\newcommand{\bcoro}{\begin{coro}}
\newcommand{\ecoro}{\end{coro}}
\newcommand{\dem}{\noindent{\bf Proof. }}

\usepackage{mathrsfs}
\renewcommand\mathcal{\mathscr}

\newcommand{\E}{{\cal E}}

\newcommand{\G}{{\cal G}}

\newcommand{\OOOO}{{\cal O}}
\renewcommand{\P}{{\cal P}}

\newcommand{\R}{{\cal R}}

\newcommand{\maths}[1]{{\mathbb #1}}  

\newcommand{\RR}{\maths{R}}
\newcommand{\NN}{\maths{N}}
\newcommand{\CC}{\maths{C}}

\newcommand{\FF}{\maths{F}}

\newcommand{\ZZ}{\maths{Z}}
\newcommand{\PP}{\maths{P}}

\newcommand{\LL}{\maths{L}}
\newcommand{\TT}{\maths{T}}

\newcommand{\uA}{\underline{A}}
\newcommand{\uD}{\underline{D}}
\newcommand{\uG}{\underline{G}}
\newcommand{\uH}{\underline{H}}
\newcommand{\uT}{\underline{T}}
\newcommand{\uP}{\underline{P}}
\newcommand{\uU}{\underline{U}}
\newcommand{\uZ}{\underline{Z}}

\newcommand{\ra}{\rightarrow}
\newcommand{\bs}{\backslash}
\newcommand{\ov}[1]{\overline{#1}} 
\newcommand{\wt}[1]{{\widetilde{#1}}}
\newcommand{\wh}[1]{{\widehat{#1}}}
\newcommand{\ga}{\gamma}
\newcommand{\Ga}{\Gamma}

\newcommand{\cqfd}{\hfill$\Box$}

\newcommand{\card}{{\operatorname{Card}}}

\renewcommand{\Im}{{\operatorname{Im}}}

\newcommand{\id}{\operatorname{id}}
\newcommand{\pr}{\operatorname{pr}}

\newcommand{\SL}{\operatorname{SL}}
\newcommand{\GL}{\operatorname{GL}}

\newcommand{\SO}{\operatorname{SO}}
\newcommand{\PO}{\operatorname{PO}}
\newcommand{\SU}{\operatorname{SU}}
\newcommand{\PU}{\operatorname{PU}}
\newcommand{\PGL}{\operatorname{PGL}}

\newcommand{\Deg}{\operatorname{deg}}

\newcommand{\bigO}{\operatorname{O}}

\newcommand{\hec}{\operatorname{hec}}
\newcommand{\heit}{\operatorname{ht}}

\newcommand{\Om}{\Omega}
\newcommand{\cem}{$c$-escape of mass}

\newcommand{\ucem}{uniform $c$-escape of mass}

\newcommand{\tr}{\operatorname{tr}}

\newcommand{\Tr}{\operatorname{Tr}}

\newcounter{fig}

\title{Escape of mass in homogeneous dynamics 
\\ in positive characteristic}
\author{Alexander Kemarsky and Frédéric Paulin and Uri Shapira}

\begin{document}
\bibliographystyle{../alphanum}
\maketitle
\begin{abstract}
We show that in positive characteristic the homogeneous probability measure supported on a periodic orbit of the diagonal group in the space of 2-lattices, when varied along rays of Hecke trees, may behave in sharp contrast to the zero characteristic analogue; that is, that for a large set of rays the measures fail to converge to the uniform probability measure on the space of 2-lattices. More precisely, we prove that when the ray is rational there is uniform escape of mass, that there are uncountably many rays giving rise to escape of mass, and that there are rays along which the measures accumulate on measures which are not absolutely continuous with respect to the uniform measure on the space of 2-lattices.

\end{abstract}

\section{Introduction}
Let $\FF_q$ be a finite field of order a positive power $q$ of a prime
$p$, and let $K=\FF_q(Y)$ be the field of rational functions in one
variable $Y$ over $\FF_q$. Let $R_\infty=\FF_q[Y]$ be the ring of
polynomials in $Y$ over $\FF_q$, let $K_\infty=\FF_q((Y^{-1}))$ be the
field of formal Laurent series in $Y^{-1}$ over $\FF_q$ and let
$X_\infty= \PGL_2(K_\infty)/\PGL_2(R_\infty)$ be the space of
homothety classes of $R_\infty$-lattices in $K_\infty\times K_\infty$
(that is, of rank $2$ free $R_\infty$-submodules spanning the vector
plane $K_\infty\times K_\infty$ over $K_\infty$).  A point $x\in
X_\infty$ is called {\it $A_\infty$-periodic} if its orbit under the
diagonal subgroup $A_\infty$ of $\PGL_2(K_\infty)$ is compact. This
orbit $A_ \infty x$ then carries a unique $A_{\infty}$-invariant
probability measure, denoted by $\mu_x$. The aim of this paper is to
study the asymptotic behavior of these measures $\mu_x$ (and in
particular to prove unexpected escape of mass phenomena) as $x$ varies
in arithmetically defined subsets of $A_\infty$-periodic points. We
will give motivations for this problem in the second part of this
introduction.

Recall that for every $x_0\in X_\infty$ and every prime polynomial
$\nu$ in $R_\infty$, the {\it Hecke tree} $T_\nu(x_0)$ with root $x_0$
is the connected component of $x_0$ in the graph with vertex set
$X_\infty$, with an edge between the homothety classes of two
$R_\infty$-lattices $\Lambda$ and $\Lambda'$ when $\Lambda'\subset
\Lambda$ and $\Lambda'/\Lambda$ is isomorphic to $R_\infty/\nu
R_\infty$ as an $R_\infty$-module. The boundary at infinity $\Om$ of
$T_\nu(x_0)$ identifies with the projective line $\PP^1(K_\nu)$ over
the completion $K_\nu$ of $K$ associated to $\nu$, and a point of
$\Om$ is called {\it rational} if it belongs to $\PP^1(K)$. (Note that
the identification of $\Om$ with $\PP^1(K_\nu)$ is not canonical, but
the notion of rationality is well defined.) For every $\xi\in\Om$, let
$(x_n^\xi)_{n\in\NN}$ be the vertices along the geodesic ray (called a
{\it Hecke ray}) in $T_\nu(x_0)$ from $x_0$ to $\xi$.

In what follows we fix an $A_{\infty}$-periodic point $x_0$ in
$X_\infty$. Note that the vertices of the Hecke-tree $T_\nu(x_0)$ then
also have periodic $A_\infty$-orbits.  Our aim is to understand the
possible sets $\Theta_\xi$
of weak-star accumulation points of the sequences of measures
$(\mu_{x^\xi_n})_{n\in\NN}$
on $X_\infty$ associated to the vertices of the Hecke ray with
endpoint $\xi$, when $\xi$ varies in $\Omega$. For all $\xi\in\Om$ and
$c>0$, we say that
\begin{enumerate}
\item[$\bullet$] $\xi$ has {\it \cem} if there exists
  $\theta\in\Theta_\xi$ with $\theta(X_\infty)\le 1-c$,
\item[$\bullet$] $\xi$ has {\it \ucem} if for every $\theta\in
  \Theta_\xi$ we have $\theta(X_\infty)\le 1-c$.
\end{enumerate}
%

Here is a summary of our results. 

\btheo\label{efrr}
There exists $c>0$ such that any rational  $\xi\in \Om$ has \ucem.
\etheo

The following result also exhibits full espace of mass phenomena along
Hecke rays.

\btheo\label{theo:fullloss}
There exists $(p,\nu,x_0)$ such that for every rational 
$\xi\in \Om$, the zero measure belongs to $\Theta_\xi$.
\etheo

The key approach to these results (proved in Subsection
\ref{subsec:ratheckeray}) is to use the geodesic flow on the quotient
of the Bruhat-Tits tree of $(\PGL_2,K_\infty)$ (see for instance
\cite{Serre83} and Subsection \ref{subsec:BruhatTitstrees}) by the
lattice $\PGL_2(R_\infty)$. 

Theorem~\ref{efrr} proves an escape of mass
phenomenon along only countably many Hecke rays. Using the remarkable
fact that the above constant $c$ is independent of the rational Hecke
ray, we can strengthen this in the next result (see Subsection
\ref{subsec:uncountableeom}).

\btheo\label{u.cthm} There exists $c>0$ such that the set of
$\xi\in\Om$ having \cem\ is uncountable.  
\etheo

As guided by the analogy with $\PGL_2(\RR)/\PGL_2(\ZZ)$ (see below),
we could still wonder if the part of the measure which does not go to
infinity still equidistributes in $X_\infty$, that is, converges to a
measure proportional to the homogeneous measure on $X_\infty$ under
$\PGL_2(K_\infty)$. The next result proves that this is also not
always the case.

\btheo\label{nac} There exists $c'>0$ such that for every
$A_\infty$-periodic point $x\in X_\infty$, there exist $\xi\in\Om$ and
$\theta\in\Theta_\xi$ such that $c'\mu_x\leq\theta$. In particular,
$\theta$ is not absolutely continuous with respect to the homogeneous
measure on $X_\infty$.  
\etheo

We give explicit constants $c,c'$ in the above statements. We will
actually prove a stronger result, Theorem \ref{th:general} in
Subsection \ref{subsec:exotic}, which mixes the behaviors in Theorems
\ref{u.cthm} and \ref{nac}.  For this, the main tool (proved in
Subsection \ref{subsec:effequisectorsphere}) is an effective
equidistribution result of sectors of Hecke spheres in positive
characteristic, which we prove using the known exponential decay of
matrix coefficients, see for instance \cite{AthGhoPra12}. We refer for
instance to the works of Dani-Margulis \cite{DanMar93},
Clozel-Oh-Ullmo \cite{CloOhUll01}, Clozel-Ullmo \cite{CloUll04},
Eskin-Oh \cite{EskOh06b}, Benoist-Oh \cite{BenOh07} for
equidistribution results of Hecke spheresxs in zero characteristic.

As we shall see in the main body of the text, Theorems \ref{efrr},
\ref{u.cthm} and \ref{nac} are valid upon replacing $K$ by any global
function field (see also Remark \ref{rem:extension} for further
extensions). In this more general case, there are several (albeit
finitely many) ways to go to infinity in $X_\infty$, and we will give
more precise results towards which cusp of $X_\infty$ the escape of
mass occurs.

\medskip 
Considering the usual analogy between function fields and number
fields, the above results should be compared with the zero
characteristic analogue, in which the behaviour is in sharp
contrast. When $R_\infty$ is replaced by $\ZZ$, $K_\infty$ by $\RR$,
and $\nu$ by an integer prime, Aka and Shapira have proved in
\cite[Theo.~4.8]{AkaSha15} that the periodic measures $\mu_{x_n^\xi}$
along (virtually any) ray in the corresponding Hecke tree
equidistribute towards the homogeneous measure in the moduli space
$\PGL_2(\RR)/\PGL_2(\ZZ)$. This result was what motivated this work,
which turned out to have a surprisingly different outcome.

The underlying phenomenon which changes drastically when passing from
zero to positive characteristics is as follows. While in zero
characteristics the size of the orbit $A_\infty x_n^\xi$ is
exponential in $n$, in positive characteristics, it is linear in $n$
due to the presence of the Frobenius automorphism (see
Theorem~\ref{theo:lineargrowpspherorbit}).  When this is combined with
the fact that rational rays diverge in a linear speed, we get the
results regarding the escape of mass.

Although the rigidity displayed in zero characteristics completely
breaks down, as demonstrated by the above results, we still believe
that the following conjecture holds.  It implies in particular, that
the set of rays having uniform escape of mass (such as the rational rays)
is a null set.

\begin{conj}\label{conj:rigidity}
  For almost any $\xi\in\Om$ (with respect to the natural probability
  measure), the averages $\frac{1}{N+1}\sum_{n=0}^N\mu_{x_n^\xi}$
  converge to the homogeneous probability measure on $X_\infty$.
\end{conj} 

Conjecture~\ref{conj:rigidity} reflects our belief that the behaviour
along rational rays is far from generic. In fact, after some computer
experiments, we suggest the following.

\begin{conj}\label{conj:full escape}
For any rational $\xi\in\Om$, $\mu_{x_n^\xi}$ converges to the 
zero measure.
\end{conj}

This work raises many other natural questions which we plan on
studying in subsequent works.  A few examples are: Is the rationality
of the ray characterized by an uniform (or full) escape of mass? Can
we find irrational rays exhibiting an escape of mass in average? Do we
have a criterion for the convergence (or convergence in average)
towards (a multiple of) the homogeneous measure ? What is the
Hausdorff dimension of the set of $\xi$ for which Theorem \ref{u.cthm}
holds?

\medskip 
Although we are working with the dynamics of rank $1$ torus,
it is interesting to compare our results with the huge corpus of works
in dynamics on noncompact spaces, in particular locally homogeneous
ones or moduli spaces, precisely devoted to prove that there is no
escape of mass to infinity for nice sequences of probability measures
on these spaces. This is in particular the case in homogeneous
dynamics -- with real Lie groups, thereby in zero characteristic --
(see for instance \cite{EskMar04,BenQui12}) or in Teichmüller dynamics
(see for instance \cite{EskMir11,Hamenstadt10a}).

Note that an escape of mass for the diagonal group is not a feature
appearing only in positive characteristics. Over the reals, there are
examples of escape of mass for the diagonal flow: for example, in
\cite{Sarnak07} the author constructs a sequence of closed geodesics
on the modular surface which converge to the zero measure (see also
\cite{Shapira16} for similar examples in higher dimensions). We stress
though that these examples do not share the arithmetic relation
between the measures along the sequence which is present in our
results. Indeed, due to the results in \cite{AkaSha15}, such an
arithmetic relation cannot coexist with an escape of mass over the
reals.

\medskip 
As another motivation for studying the limiting behaviour of
$\mu_{x_n^\xi}$ (also originating from the analogy with
\cite{AkaSha15}), let us give a relation with the distribution
properties of the periods of the continued fraction expansion of
certain sequences of quadratic irrationals. We refer for instance to
the surveys \cite{Lasjaunias00,Schmidt00} for background. We denote by
$O_\infty=\FF_q[[Y^{-1}]]$ the local ring of $K_\infty$ (consisting of
power series in $Y^{-1}$ over $\FF_q$). Any element $f\in K_\infty$
may be uniquely written $f=[f]+\{f\}$ with $[f]$ in the polynomial
ring $R_\infty=\FF_q[Y]$ and $\{f\}\in Y^{-1} O_\infty$. The {\it
  Artin map} $\Psi:Y^{-1} O_\infty-\{0\}\ra Y^{-1}O_\infty$ is defined
by $f\mapsto \big\{\frac{1}{f}\big\}$. Any $f\in K_\infty$ irrational
(not in $K=\FF_q(Y)\,$) has a unique continued fraction expansion
$$
f=a_0 +
\cfrac{1}{a_1+\cfrac{1}{a_2+ \cfrac{1}{a_3+\cdots}}}\;,
$$
with $a_0=[f]\in R_\infty$ and $a_n= \big[\frac{1}{\Psi^{n-1}(f-a_0)}
\big]$ a non constant polynomial, for $n\geq 1$. Let $QI=\{f\in
K_\infty\;:\;[K(f):K]=2\}$ be the set of quadratic irrationals over
$K$ in $K_\infty$. Assume for simplicity that the characteristic $p$
is different from $2$, and denote by $f^\sigma\neq f$ the Galois
conjugate of $f\in QI$ over $K$.  Given an irrational $f\in Y^{-1}
O_\infty$, we have $f\in QI$ if and only if the continued fraction
expansion is eventually periodic. We then denote by $\nu_f$ the
uniform probability on the periodic part of the orbit of $f$ under
$\Psi$, by $g_f=
\begin{bmatrix} f & f^\sigma \\ 1 & 1\end{bmatrix}\in
\PGL_2(K_\infty)$, and by $x_f=g_f\PGL_2(R_\infty)\in X_\infty$. It is
then easy to prove that $x_f$ is $A_\infty$-periodic. Using the main
results of \cite{BroPau07JLMS}, we may construct a natural map from (a
full-measure subset of) $X_\infty$ onto (a full-measure subset of)
$Y^{-1} O_\infty$, sending the (normalized) homogeneous measure
$m_\infty$ on $X_\infty$ to the (normalized) Haar measure on $Y^{-1}
O_\infty$, $A_\infty$-orbits in $X_\infty$ to $\Psi$-orbits in $Y^{-1}
O_\infty$, and more precisely the $A_\infty$-invariant probability
measure $\mu_{x_f}$ to the equiprobability $\nu_f$ for every quadratic
irrational $f$ in $Y^{-1} O_\infty$. Hence the distribution properties
of the periods of the continued fraction expansions of quadratic
irrationals are related to the distribution properties of the
$A_\infty$-orbits in $X_\infty$.  

\medskip\noindent{\small {\it Acknowledgements: } We thank the
  hospitality of the Institut Henri Poincaré in early 2014 where part
  of this work was done. This work was supported by the NSF Grant no
  093207800, while the last two authors were in residence at the MSRI,
  Berkeley CA, during the Spring 2015 semester. We thank J.-F.~Quint
  (for his help for the proof of Proposition
  \ref{prop:unipradabetransitif}), Y.~Benoist, L.~Clozel,
  G.~Chenevier, M.~Einsiedler, and E.~Lindenstrauss for discussions on this paper. U.~S.~acknowledges the
  support of ISF grant 357/13.}

\section{Global function fields and Bruhat-Tits trees}
\label{sec:background} 

This section introduces the notation and preliminary results used in
this paper.  We refer the reader to the following commutative diagram
for a global view of this notation.

\bigskip
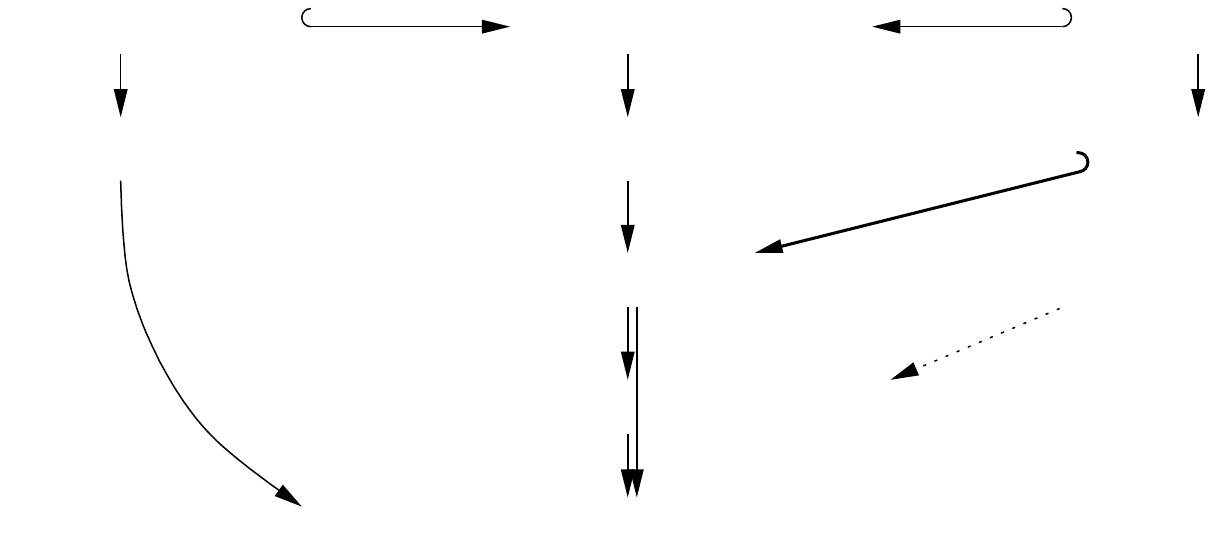

\subsection{Global function fields}
\label{subsec:function fields}

We refer for instance to \cite{Rosen02,Serre79} for the content of
this subsection. 

Let $\FF_q$ be a finite field with $q$ elements, where $q$ is a
positive power of a prime $p$. Let $K$ be a {\it global function
  field} over $\FF_q$, that is, the function field of a geometrically
connected smooth projective curve $\bf C$ over $\FF_q$, or
equivalently an extension of $\FF_q$ of transcendance degree $1$, in
which $\FF_q$ is algebraically closed. The set $\P$ of {\it primes} of
$K$ is the set of closed points of $\bf C$, or equivalently the set of
discrete valuations of $K$, trivial on $\FF_q^\times$, with value
group exactly $\ZZ$. We fix an element in $\P$ that we denote by
$\infty$, and we denote by $\P_f$ the set $\P-\{\infty\}$.

For every $\omega\in\P$, we denote by $R_\omega$ the affine algebra of
the affine curve ${\bf C}-\{\omega\}$ (which is a Dedekind ring), by
$v_\omega$ the discrete valuation of $K$ associated to $\omega$ (with
the usual convention that $v_\omega(0)=+\infty$), by $K_\omega$ the
associated completion of $K$ (and again by $v_\omega$ the extension of
$v_\omega$ to $K_\omega$), by $O_\omega$ its local ring, by
$\pi_\omega$ a uniformizer of $O_\omega$, by $k_\omega$ its
residual field (that we identify with its canonical lift in
$O_\omega$), and by $\Deg(\omega)$ the degree of $k_\omega$ over
$\FF_q$.  We assume, as we may using for instance the Riemann-Roch
theorem,
%
%
that $\pi_\nu$ belongs to $R_\infty$ if $\nu\in\P_f$.  Note that
$R_\infty\subset O_\nu$ if $\nu\in\P_f$ (since an element in
$R_\infty$ has no pole at the closed point $\nu\neq \infty$ of ${\bf
  C}$), and that $R_\infty[\pi_\nu^{-1}]\cap O_\nu= R_\infty$.

We normalize the absolute value $|\cdot|_\omega$ associated to
$v_\omega$ by $|x|_\omega= |k_\omega|^{-v_\omega(x)}= q^{-\Deg
  \omega\; v_\omega(x)}$ for every $x\in K_\omega$. In particular, the
product formula
$$
\forall\; x\in K,\;\;\;\prod_{\omega\in\P}|x|_\omega=1
$$
holds. Note that $K_\omega$ is the field $k_\omega((\pi_\omega))$ of Laurent
series $f= \sum_{i\in\ZZ} f_i(\pi_\omega)^{i}$ in the variable $\pi_\omega$
over $k_\omega$, where $f_i\in k_\omega$ is zero for $i\in\ZZ$ small
enough. We have
%
%
$$
|f|_\omega=|k_\omega|^{-\sup\{j\in\ZZ\;:\;\forall\;i<j,\;f_i=0\}}\;,
$$
and $O_\omega=k_\omega[[\pi_\omega]]$ is the local ring of power series
$f=\sum_{i\in\NN} f_i(\pi_\omega)^{i}$ (where $f_i\in k_\omega$) in the
variable $\pi_\omega$ over $k_\omega$. 

For every finite extension $\wt K_\omega$ of $K_\omega$, we denote
again by $v_\omega$ the unique extension of $v_\omega$ to a valuation
on $\wt K_\omega$, and by $e(\wt K_\omega, K_\omega) = [v_\omega({\wt
  K_\omega}^\times): v_\omega({K_\omega}^\times)]$ its ramification
index (see for instance \cite[\S 2]{Serre79}).

For instance, if $\bf C$ is the projective line $\PP^1$ and if
$\infty=[1:0]$ is its usual point at infinity, then $K=\FF_q(Y)$,
$\pi_\infty=Y^{-1}$, $K_\infty= \FF_q((Y^{-1}))$, $O_\infty=
\FF_q[[Y^{-1}]]$, $k_\infty=\FF_q$, $R_\infty=\FF_q[Y]$ and the
uniformizers $\pi_\nu$ for $\nu\in\P_f$ may be taken to be the monic
prime polynomials in $R_\infty$, with $\Deg \nu$ the degree of the
polynomial $\pi_\nu$. This is the example considered in the
introduction.

\subsection{Generalisation to rank-one semi-simple groups}
\label{subsec:algebraic groups}

The aim of this subsection is to explain to which $K_\infty$-rank-one
groups the tools introduced in this paper are applying besides
$\PGL_2$. But for the readability, we will restrict to this last case
at the end of this subsection, giving the group-theoretic notation we
are going to use. We refer for instance to \cite{Tits66,Tits79} for
the already known content of this subsection.

Let $\uG$ be a connected semi-simple linear algebraic group defined
over $K$, with $K_\infty$-rank one. We fix an embedding $\uG\ra \GL_N$
for some $N\in\NN$. The example considered in the introduction is
$\uG=\PGL_2$ (which is adjoint and absolutely simple). 

For every $\omega\in\P$ and every algebraic subgroup $\uH$ of $\uG$
defined over $K_\omega$ (for instance if $\uH$ is defined over $K$), we
set $H_\omega= \uH(K_\omega)$, which is a non-Archimedian Lie group.

For every $\omega\in\P$, we define $\Ga_\omega=\uG(R_\omega)$, which is a
lattice in the locally compact group $G_\omega$. For instance, when
$\uG=\PGL_2$ and ${\bf C}=\PP^1$, the lattice $\Ga_\infty$ is called
{\it Nagao's lattice} \cite{Nagao59} (or Weil's modular group
\cite{Weil70}).

For every $\omega\in\P$, we denote by $X_\omega$ the totally
disconnected locally compact space $\Ga_\omega\bs G_\omega$
(contrarily to the introduction, we consider the left quotient, as it
makes the connection with Bruhat-Tits theory easier). As we want to
study phenomena of escape of mass at infinity for measures on
$X_\infty$, we require that $X_\infty$ is not compact.  For instance,
when $\uG=\PGL_2$, the space $X_\infty$ is non compact, and identifies
by $\Ga_\infty g\mapsto g^{-1} [R_\infty\times R_\infty]$ with the
space of homothety classes $[\Lambda]$ under $K_\infty^\times$ of
$R_\infty$-lattices $\Lambda$ in $K_\infty\times K_\infty$.

Given $\nu\in\P_f$, let $S=\{\infty,\nu\}$ and let$\Ga_S$ be the {\it
  $S$-arithmetic group} $\uG(R_\infty[{\pi_\nu}^{-1}])$, which embeds
diagonally in the locally compact group $G_S=G_\infty\times G_\nu$ as
a lattice, and let $X_S=\Ga_S\bs G_S$. We identify $G_\infty$ and
$G_\nu$, hence any subgroup of them, with their images in $G_S$ by the
maps $x\mapsto (x,e)$ and $y\mapsto (e,y)$.

Note that when $\nu\in \P_f$, the $K_\nu$-rank of $\uG$ may be $1$ (as
in the case $\uG=\PGL_2$) or not. For instance, let $D$ be a
(finite dimensional) central simple algebra over $K$ which is {\it
  ramified at $\infty$} (that is, $D_\infty=D\otimes_K
K_\infty)$ is a division algebra). Then the algebraic group $\uG$ with
$\uG(L)=\PGL_2(D\otimes_KL)$ for every $K$-algebra $L$ is an
(adjoint absolutely quasi-simple) connected semi-simple linear
algebraic group defined over $K$, with $K_\infty$-rank one. For all
$\nu\in \P_f$, the group $\uG$ has $K_{\nu}$-rank $1$ if and only if
$D$ ramifies at $\nu$ (that is, when $D_\nu= D \otimes_K
K_\nu$ is a division algebra).

\medskip The next two results are not necessary for the main results
of the paper, but they will be used to explain the restrictions on the
considered algebraic groups.  The first one follows from a well-known
argument of weak approximation.

\blemm \label{lem:bisplittoriexist} Let $\nu\in \P_f$, if the
$K_\nu$-rank of $\uG$ is $1$, then there exist tori $\uA$ in $\uG$
defined over $K$, which splits over both $K_\infty$ and $K_\nu$ (hence
is a maximal $K_\infty$-split and $K_\nu$-split torus).  
\elemm

\dem By \cite[Theo.~2]{PraRap02} applied to the semisimple connected
algebraic group $\uG$ defined over the infinite field $K$, there
exists $m\in\NN-\{0\}$ such that the closure of the image of the
diagonal embedding of $\uG(K)$ in $G_\infty\times G_\nu$ contains the
subgroup of $G_\infty\times G_\nu$ generated by $m$-th powers. Let
$\ga_\infty$ and $\ga_\nu$ be nontrivial elements in $\uG(K)$ which
split over $K_\infty$ and $K_\nu$ respectively. There hence exists an
element in $\uG(K)$ arbitrarily close to both $\ga_\infty^m$ and
$\ga_\nu^m$, which therefore splits simultaneously over $K_\infty$ and
$K_\nu$.  \cqfd

\bprop\label{prop:unipradabetransitif} Let $\uH$ be an adjoint,
absolutely quasi-simple, connected, semi-simple algebraic group over a
local field $F$ of $F$-rank one. Let $\uT$ be a maximal $F$-split
torus, $\uZ$ its centralizer, $\uP$ a minimal parabolic subgroup of
$\uH$ over $F$, and $\uU$ its unipotent radical. If $\uH$ is
isomorphic over $F$ to the algebraic group $L\mapsto
\PGL_2(D\otimes_FL)$ for every $F$-algebra $L$, where $D$ is a central
division algebra over $F$, then $\uZ(F)$ acts transitively on
$\uU(F)-\{0\}$ by conjugation. If $\uH$ is isomorphic over $F$ to the
algebraic group $L\mapsto \PU_{1,1}(D\otimes_FL)$ for every
$F$-algebra $L$, where $D$ is a quaternion division algebra over $F$,
then $\uZ(F)$ acts non transitively with finitely many orbits on
$\uU(F)-\{0\}$ by conjugation. Otherwise, $\uZ(F)$ acts with
infinitely many orbits on $\uU(F)-\{0\}$ by conjugation.  \eprop

In particular, by the classification theorem \cite{Tits66}, if
furthermore $F=K_\nu$ for some $\nu\in\P$ and $\uH$ is defined and
isotropic over $K$, then $\uZ(F)$ acts transitively on $\uU(F)-\{0\}$
by conjugation if and only if $\uH$ is isomorphic over $K$ to
$\PGL_2(\uD)$ where $\uD$ is a central division algebra over $K$, and
$\uZ(F)$ acts non transitively with finitely many orbits on
$\uU(F)-\{0\}$ by conjugation if and only if $\uH$ is isomorphic over
$K$ to $\PU_{1,1}(\uD)$ where $\uD$ is a quaternion division algebra
over $K$.

\medskip 
\dem Let $H=\uH(F)$, $T=\uT(F)$, $Z=\uZ(F)$ and $U=\uU(F)$.
Let us denote by $\begin{bmatrix} a_{ij}\end{bmatrix}$ the image in
$\PGL_2$ of a matrix $\begin{pmatrix} a_{ij}\end{pmatrix}$ in $\GL_2$.

If $\uU$ is non abelien (or equivalently if the (relative) root system
of $\uH$ is not reduced), it is easy to see that the action of
$\uZ(F)$ on $\uU(F)-\{0\}$ by conjugation has infinitely many orbits.
Conversely, assume that $\uU$ is abelien.  When $F=\CC$ (then
$H=\PGL_2(\CC)\,$) or  $F=\RR$ (then $H=\PO(1,n)$),  the
action of $\uZ(F)$ on $\uU(F)-\{0\}$ by conjugation is
transitive. Hence assume that $F$ is non Archimedian.  By the
classification theorem \cite{Tits79}, up to isomorphism, $H$ is either
$\PGL_2(D)$ for a central division algebra $D$ over $F$, or
$\PU_{1,1}(D)$ for a quaternion division algebra $D$ over $F$ and the
Hermitian form $h(z_1,z_2)= \overline{z_1}z_2+\overline{z_2}z_1$.

\medskip
In the first of the above two cases, we may take
$$
T= \Big\{\begin{bmatrix} a & 0 \\
  0 & d\end{bmatrix}\;:\;a,d\in F^\times\Big\} ,\;
Z= \Big\{\begin{bmatrix} a & 0 \\
  0 & d\end{bmatrix}\;:\;a,d\in D^\times\Big\} ,\;
U= \Big\{\begin{bmatrix} 1 & b \\
  0 & 1\end{bmatrix}\;:\;b\in D\Big\}\;.
$$
The transitivity of the action by conjugation of $Z$ on $U-\{0\}$
follows hence from the transitivity of the action of $D^\times\times
D^\times$ on $D-\{0\}$ by $(a,d)\cdot b= abd^{-1}$, which is
immediate.  

\medskip In the second case, we denote by $z\mapsto \overline{z}$ the
canonical involution in the quaternion division algebra $D$ over $F$,
by $N:x\mapsto x\overline{x}$ and $\Tr:x\mapsto x+\overline{x}$ its
(reduced) norm and trace, and by $(1,i,j,k)$ a standard basis of $D$
over $F$. Recall that $F^\times/(F^\times)^2$ is finite and non
trivial. Indeed, this group is isomorphic to $(\ZZ/2\ZZ)\times
\big(f^\times/ (f^\times)^2 \big)$ where $f$ is the (finite) residue
field of $F$, since if $\OOOO_F$ is the local ring and $\pi_F$ is a
uniformizer in $F$, the map $(n,x)\mapsto \pi_F^nx$ from $\ZZ\times
\OOOO_F^\times$ to $F^\times$ is an isomorphism.

Let $\Im\;D=\{x\in D\;:\;\Tr(x)=0\}$ be the $K$-vector space of purely
imaginary elements of $D$, endowed with the action of the orthogonal
group $\operatorname{O}(N_{\mid\Im\;D})$ of the restriction to
$\Im\;D$ of the norm. Since $F^\times/(F^\times)^2$ is finite and
$N(F^\times) =(F^\times)^2$, there exists a finite subset $A$ of
$F^\times$ such that every line in $\Im\;D$ contains a vector whose
norm lies in $A$. By Witt's theorem, the group
$\operatorname{O}(N_{\mid\Im\;D})$ hence acts with finitely many
orbits on the lines of $\Im\;D$.

The group $\SL_2(D)$ acts by $g\cdot M={\,}^t\overline{g}\,Mg$ on the
$6$-dimensional $F$-vector space $E=\big\{M=\begin{pmatrix} a & b\\
  \overline{b} & d\end{pmatrix} \;:\; a,d\in F, b\in D\big\}$, by
preserving the Dieudonné determinant $\det M= ac-N(b)$, which is a
quadratic form $Q$ on $E$.  Let $M_0=\begin{pmatrix} 0 & 1\\
  1 & 0\end{pmatrix}$, which is a $Q$-anisotropic element of $E$. The
group $\SU_{1,1}(D)$ is the stabilizer of $M_0$ in $\SL_2(D)$ for the
above action. Let $M_0^{\perp}$ be the $5$-dimensional orthogonal of
$M_0$ in $E$ for $Q$, which is invariant under $\SU_{1,1}(D)$, and
note that the restriction $Q_{\mid M_0^{\perp}}$ is non degenerate.
We consider the basis
$$
\Big(e_1=\begin{pmatrix} 1 & 0\\
  0 & 0\end{pmatrix}, e_2=\begin{pmatrix} 0 & i\\
  \overline{i} & 0\end{pmatrix}, e_3=\begin{pmatrix} 0 & j\\
  \overline{j} & 0\end{pmatrix}, e_4=\begin{pmatrix} 0 & k\\
  \overline{k} & 0\end{pmatrix}, e_5=\begin{pmatrix} 0 & 0\\
  0 & 1\end{pmatrix}\Big)
$$
of $M_0^{\perp}$, and we sometimes write matrices by blocs in the
decomposition $(e_1, (e_2,e_3,e_4), e_5)$. The group $T=
\Big\{\begin{bmatrix} \lambda & 0 & 0\\ 0 & \id & 0\\ 0 & 0 &
  \lambda^{-1}\end{bmatrix}\;: \;\lambda\in F^*\Big\}$ is a 
maximal $F$-split torus in $\PO( Q_{\mid M_0^{\perp}})$, whose
centralizer $Z$ contains $Z'=\Big\{\begin{bmatrix} 
1 & 0 & 0\\ 0 & A &  0\\ 0 & 0 & 1\end{bmatrix}\;:
\;A\in\operatorname{O}(N_{\mid\Im\;D})\Big\}$.  The
projective upper triangular subgroup $P$ of $\PO( Q_{\mid
  M_0^{\perp}})$ is a minimal $F$-parabolic subgroup of $\PO( Q_{\mid
  M_0^{\perp}})$, whose unipotent radical is, by an easy computation,
$$
U=\Big\{\begin{bmatrix} 
1 & x N(i) & y N(j) & z N(k) & N(xi+yj+zk)\\ 
0 & 1 & 0 & 0 & x\\0 & 0 & 1 & 0 & y\\0 & 0 & 0 & 1 & z\\ 
0 & 0 &0 & 0 & 1\end{bmatrix}\;:\;x,y,z\in F\Big\}\;.
$$
The action of $Z'$ on $U$ by conjugation thus identifies with the
linear action of $\operatorname{O}(N_{\mid\Im\;D})$ on $\Im\;D$.  The
natural map $\SU_{1,1}(D)\ra \SO( Q_{\mid M_0^{\perp}})$ induced by
the action of $\SU_{1,1}(D)$ on $M_0^{\perp}$ is an isogeny, by the
semi-simplicity of $\SU_{1,1}(D)$. Hence the adjoint groups
$\PU_{1,1}(D)$ and $\PO( Q_{\mid M_0^{\perp}})$ are isomorphic. 

For every $\lambda\in F$, the action of $\begin{bmatrix} \lambda & 0 &
  0\\ 0 & \id & 0\\ 0 & 0 & \lambda^{-1}\end{bmatrix}\in T$ on each
line in $U$ is the multiplication by $\lambda^2$. Since $F^\times/
(F^\times)^2$ is finite (of cardinality at least $2$), the action of
$T$ on each line in $U$ has finitely many orbits (and at least
two). The fact that the action of $\uZ(F)$ on $\uU(F)-\{0\}$ by
conjugation has finitely many orbits (and at least two) hence follows
from the fact that the action of $\operatorname{O} (N_{\mid\Im\;D})$
on the lines of $\Im\;D$ has finitely many orbits.  
\cqfd

\bigskip From now on in this paper, we fix $\nu\in\P_f$ and we denote
by $\uG=\PGL_2$ the (adjoint semi-simple absolutely simple) projective
linear algebraic group over $K$ in dimension $2$, so that
$\Ga_\infty=\uG(R_\infty)= \PGL_2(R_\infty)$ is a nonuniform lattice
in $G_\infty= \PGL_2(K_\infty)$. Whenever necessary, we embed $\PGL_2$
in $\GL_3$ by the adjoint representation on the vector space of
traceless $2$-by-$2$ matrices.

Let $\uA$ be the {\it diagonal subgroup} of $\uG$, that is, the
algebraic subgroup of $\uG$ consisting in the elements represented by
diagonal matrices, which is a (split) maximal torus  of $\uG$
defined over $K$.

Let $S=\{\infty,\nu\}$, so that the $S$-arithmetic group $\Ga_S=
\uG(R_\infty[{\pi_\nu}^{-1}])$ is a nonuniform lattice in
$G_S=G_\infty\times G_\nu$. Note that $\Ga_S\cap\uG(O_\nu)=\Ga_\infty$
since $R_\infty[\pi_\nu^{-1}]\cap O_\nu= R_\infty$.

For every $\omega\in S$, we denote

$\bullet$~ by $\begin{bmatrix} a & b \\ c & d \end{bmatrix}$ the image
in $G_{\omega}= \PGL_2(K_{\omega})$ of $\begin{pmatrix} a & b \\ c &
  d \end{pmatrix}\in \GL_2(K_{\omega})$,

$\bullet$~ again by $v_\omega$ the map from the abelian group
$A_{\omega}= \uA(K_{\omega})=\Big\{\begin{bmatrix} a & 0 \\ 0 &
  d \end{bmatrix}\;:\;a,d\in K_{\omega}\Big\}$ to $\ZZ$ defined by
$\begin{bmatrix} a & 0 \\ 0 & d \end{bmatrix}\mapsto v_\omega(d/a)$,
which is a group epimorphism with compact-open kernel
$\uA(O_\omega)=\uG(O_\omega)\cap A_\omega$,

$\bullet$~ by $\alpha_\omega:K_\omega^\times\ra A_\omega$ the group
isomorphism $t\mapsto \begin{bmatrix} 1 & 0 \\ 0 & t \end{bmatrix}$
(whose inverse is the positive root of the torus $\uA$ over
$K_\omega$), so that $v_\omega(\alpha_\omega(t))=v_\omega(t)$,
and by $a_\omega=\alpha_\omega(\pi_\omega)=\begin{bmatrix} 1 & 0 \\
  0 & \pi_\omega \end{bmatrix}$, so that $v_\omega(a_\omega)=1$.

\brema\label{rem:extension}{\rm  An appropriate version of this paper
(including loss of mass phenomena of the homogeneous probability
measures on the periodic orbits of the points along appropriate rays
of the Hecke tree of any given periodic point of $X_\infty$) is valid
when we replace $\uG$ by the linear algebraic group over $K$ defined

\begin{enumerate}
\item[$\bullet$] either by $\uG(L)=\PGL_2(D\otimes_KL)$ for every
$K$-algebra $L$, where $D$ is a (finite dimensional) central division
algebra over $K$ which ramifies at the places $\infty$ and $\nu$, and
we endow the algebraic group $\uG$ with a $R_\infty$-structure such
that $\uG(R_\infty)= \PGL_2(\R_\infty)$ where $\R_\infty$ is a
$R_\infty$-order in $D$ (see \cite{Reiner75} for any information on
orders),

\item[$\bullet$] or by $\uG(L)=\PU_{1,1}(D\otimes_KL)$ for every $K$-algebra
$L$, where $D$ is a quaternion algebra over $K$ (and the underlying
Hermitian form is $(z_1,z_2)\mapsto
\overline{z_1}z_2+\overline{z_2}z_1$),
\end{enumerate}
\noindent allowing, thanks to the transitivity properties described in
Proposition \ref{prop:unipradabetransitif}, to prove a modified
version of Theorem \ref{theo:lineargrowpspherorbit}, when we replace
$\Ga_\infty$ by a congruence subgroup and when we replace $\uA$ by any
torus over $K$ in $\uG$ which splits over both $K_\infty$ and $K_\nu$
(which exists by Lemma \ref{lem:bisplittoriexist}). But for the sake
of simplicity, we stick to the above choice of
$(\nu,\uG,\uA,\Ga_\infty)$.}
\erema

\subsection{Bruhat-Tits trees}
\label{subsec:BruhatTitstrees}

Let $(K,\nu,\uG,\uA)$ and the associated notation be as in Subsection
\ref{subsec:algebraic groups} before Remark \ref{rem:extension}.

\medskip
\noindent{\bf Trees. }
Let $T$ be a locally finite tree. Its set of vertices $VT$ is endowed
with the maximal distance for which two adjacent distinct vertices are
at distance $1$. A {\it geodesic ray} or {\it line} in $T$ is an
isometric map from $\NN$ or $\ZZ$ to its set of vertices. The set of
geodesic lines of $T$, endowed with the compact-open topology, is denoted
by $\G T$.

An {\it end} of $T$ is an equivalence class of geodesic rays, when two
geodesic rays are equivalent if the intersection of their images is
the image of a geodesic ray. The set of ends of $T$, endowed with the
(compact, totally disconnected) quotient topology of the compact-open
topology, is denoted by $\partial_\infty T$, and called the {\it
  boundary at infinity} of $T$.

The {\it translation length} of an isometry $\ga$ of $T$ is
$$
\ell_T(\ga)=\min_{x\in VT} d(x,\ga x)\;.
$$
It is invariant under conjugation of $\ga$ in the isometry group of
$T$. We will say that $\ga$ is {\it loxodromic} if $\ell_T(\ga)>0$, in
which case there exists a unique image of a geodesic line in $T$
on which $\ga$ translates a distance $\ell_T(\ga)$, called the {\it
  translation axis} of $\ga$.

The {\it geodesic flow} (with discrete times) $(\phi_m)_{m\in\ZZ}$ on
$T$ is the right action $(\G T\times \ZZ)\ra \G T$ of $\ZZ$ on $\G T$
by translations at the source, defined by 
$$
(\ell,m)\mapsto \{\phi_m\ell:n\mapsto \ell(n+m)\}
$$ 
for all $m\in\ZZ$ and $\ell:\ZZ\ra VT$ in $\G T$.  Given a group $\Ga$ of
automorphisms of $T$, the geodesic flow on $T$ induces a right action
of $\ZZ$ on $\Ga\bs \G T$, also called the {\it geodesic flow} of $\Ga
\bs T$, and again denoted by $(\phi_m)_{m\in\ZZ}$.

\medskip
\noindent{\bf The tree of $\PGL_2$ over local fields. }  For
$\omega\in S=\{\infty,\nu\}$, let $\TT_\omega$ be the {\it Bruhat-Tits
  tree} of $(\uG,K_\omega)$, see for instance \cite{Tits79}. We use
its description given in \cite{Serre83}.

Recall that an {\it $O_{\omega}$-lattice} $\Lambda$ in the
$K_\omega$-vector space $K_\omega\times K_\omega$ is a rank $2$ free
$O_{\omega}$-submodule of $K_\omega\times K_\omega$, generating
$K_\omega\times K_\omega$ as a vector space. The Bruhat-Tits tree
$\TT_\omega$ is the graph whose set of vertices $V\TT_\omega$ is the
set of homothety classes (under $K_\omega^\times$) $[\Lambda]$ of
$O_{\omega}$-lattices $\Lambda$ in $K_\omega\times K_\omega$, and
whose non-oriented edges are the pairs $\{x,x'\}$ of vertices such
that there exist representatives $\Lambda$ of $x$ and $\Lambda'$ of
$x'$ such that $\Lambda \subset \Lambda'$ and $\Lambda'/\Lambda$ is
isomorphic to $O_{\omega}/{\pi_{\omega}} O_{\omega}$.  This graph is a
regular tree of degree $|\PP_1(k_{\omega})|= |k_{\omega}|+1$.

We denote by $*_\omega$ the homothety class of the
$O_{\omega}$-lattice $O_{\omega} \times O_{\omega}$ generated by the
canonical basis of $K_\omega\times K_\omega$.  The left linear action
of $\GL_2(K_\omega)$ on $K_\omega\times K_\omega$ induces a faithful,
transitive left action of $G_\omega$ on $V\TT_\omega$. The stabilizer
in $G_\omega$ of $*_\omega$ is $\uG(O_{\omega})$.  We will hence
identify $G_\omega/ \uG(O_{\omega})$ with $V\TT_\omega$ by the map
$g\,\uG(O_{\omega})\mapsto g\,*_\omega$.

We identify as usual the projective line $\PP_1(K_\omega)$ with
$K_\omega \cup \{\infty\}$ using the map $K_\omega(x,y) \mapsto
xy^{-1}$. There exists one and only one homeomorphism between the
boundary at infinity $\partial_\infty \TT_\omega$ of $\TT_\omega$ and
$\PP_1(K_\omega)$ such that the (continuous) extension to
$\partial_\infty \TT_\omega$ of the isometric action of $G_\omega$ on
$\TT_\omega$ corresponds to the projective action of $G_\omega$ on
$\PP_1(K_\omega)$. From now on, we identify $\partial_\infty
\TT_\omega$ and $\PP_1(K_\omega)$ by this homeomorphism.

The group $G_\omega$ hence acts simply transitively on the set of
ordered triples of distinct points in $\partial_\infty \TT_\omega$. In
particular, the group $G_\omega$ acts transitively on the space
$\G\TT_\omega$ of geodesic lines in $\TT_\omega$. The stabilizer under
this action of the geodesic line
$$
\ell_0:n\mapsto 
[O_\omega\times\pi_\omega^{n}O_\omega]=a_\omega^n\;*_\omega
$$ 
is the maximal compact-open subgroup $\uA(O_\omega)$ of the diagonal
group $A_\omega$. We will hence identify $G_\omega/\uA(O_\omega)$ with
$\G\TT_\omega$ by $g\uA(O_\omega)\mapsto g\,\ell_0$. Furthermore, the
stabilizer in $G_\omega$ of the ordered pair of endpoints
$(\ell_0(-\infty)=0, \ell_0(+\infty)=\infty)$ of $\ell_0$ in
$\partial_\infty \TT_\omega=\PP_1(K_\omega)$ is $\uA_\omega$.
Therefore any element $\ga_0\in G_\omega$ which is loxodromic on
$\TT_\omega$ is diagonalisable over $K_\omega$. Besides, by \cite[page
108]{Serre83}, the translation length on $\TT_\omega$ of $\ga_0= 
\begin{bmatrix} \lambda_+ & 0 \\ 0 & \lambda_- \end{bmatrix}$ is
\begin{equation}\label{eq:translatlength}
\ell_{\TT_\omega}(\ga_0)= |v_\omega(\lambda_+)-v_\omega(\lambda_-)|\;.
\end{equation}
Using the group morphism $v_\omega:A_\omega\ra\ZZ$, the action by
translations on the right of $A_\omega$ on $G_\omega/\uA(O_\omega)$
corresponds to the geodesic flow on $\G\TT_\omega$: for all $a\in A_\omega$
and $\ell\in \G\TT_\omega= G_\omega/\uA(O_\omega)$, we have 
$$
\ell\;a=\phi_{v_\omega(a)} \ell\;. 
$$

We denote by $\pi'_\infty : X_\infty= \Ga_\infty \bs G_\infty\ra
\Ga_\infty\bs \G\TT_\infty=\Ga_\infty\bs G_\infty /\uA(O_\infty)$ the
canonical projection (see the diagram at the beginning of Section
\ref{sec:background}). The previous equation proves that $\pi'_\infty$
is equivariant with respect to the morphism $v_\infty:A_\infty\ra\ZZ$,
where $A_\infty$ acts by translation on the right on $X_\infty$ and
$\ZZ$ by the (quotient) geodesic flow on $\Ga_\infty\bs \G\TT_\infty$:
for all $x\in X_\infty$ and $a\in A$,
\begin{equation}\label{eq:equivarvomega}
\pi'_\infty(x \,a)=\phi_{v_\infty(a)}\pi'_\infty(x)\;.
\end{equation} 

\bigskip
\noindent {\bf The principal bundle $\pi_\infty: X_S\ra X_\infty$. }
Since $\Ga_S$ is irreducible, the group $\Ga_\infty=\Ga_S\cap
\uG(O_\nu)$ is dense in the stabiliser $\uG(O_\nu)$ of the base point
$*_\nu$ of the Bruhat-Tits tree $\TT_\nu$.  This stabilizer
$\uG(O_\nu)$ acts transitively on the geodesic rays in $\TT_\nu$
starting from $*_\nu$. Thus $\Ga_\infty$ preserves and acts
transitively on the sphere in $\TT_\nu$ of any given radius centered
at $*_\nu$. For every $g'\in G_\nu$, there hence exists $\ga\in
\Ga_\infty$ and $n\in\NN$ such that $\ga^{-1}g' *_\nu= [O_\nu \times
\pi_\nu^{n}O_\nu]= a_\nu^{n}\; *_\nu$. Therefore
\begin{equation}\label{eq:decomradialGnu}
G_\nu= \bigcup_{n\in\NN} \Ga_\infty \;a_\nu^{n}\; \uG(O_\nu)\;.
\end{equation}  
In particular, $G_\nu=\Ga_S\; \uG(O_\nu)$.

Therefore, every element $x$ of $X_S$ may be written $\Ga_S(g,g')$
with $g\in G_\infty$ and $g'\in \uG(O_\nu)$. For all $g,h\in G_\infty$
and $g',h'\in \uG(O_\nu)$, we have $\Ga_S(g,g')=\Ga_S(h,h')$ if and
only if $gh^{-1}=g'(h')^{-1}\in\Ga_S\cap \uG(O_\nu)=\Ga_\infty$. Hence
the map $\pi_\infty:X_S\ra X_\infty$, where $\pi_\infty(x)=\Ga_\infty
g$ if $x=\Ga_S(g,g')$ with $g'\in \uG(O_\nu)$, is well defined and
continuous. The action of $\uG(O_\nu)$ by right translations on the
second factor of $G_S=G_\infty\times G_\nu$ induces an action of $\uG
(O_\nu)$ on $X_S=\Ga_S\bs G_S$, which is transitive and free on the
fibers of $\pi_\infty$. Hence $\pi_\infty: X_S\ra X_\infty$ is a
principal bundle under the group $\uG(O_\nu)$, which gives an
identification between $X_\infty=\Ga_\infty\bs G_\infty$ and
$X_S/\uG(O_\nu)= \Ga_S \bs G_S/\uG(O_\nu)$ (see the diagram at the
beginning of Section \ref{sec:background}).

\bigskip
\noindent{\bf Ends of the modular graph at the place $\infty$ 
and heights. }
The quotient graph $\Ga_\infty\bs \TT_\infty$ will be called the {\it
  modular graph at $\infty$} of $K$. By for instance \cite{Serre83},
the {\it set of cusps} $\Ga_\infty\bs\PP_1(K)$ is finite, and
$\Ga_\infty\bs \TT_\infty$ is the disjoint union of a finite connected
subgraph containing $\Ga_\infty*_\infty$ and of maximal open geodesic
rays $h_z (\,]0,+\infty[)$, for $z=\Ga_\infty\wt z\in \Ga_\infty\bs
\PP_1(K)$, where $h_z$ (called a {\it cuspidal ray}) is the image by
the canonical projection $\TT_\infty\ra \Ga_\infty \bs \TT_\infty$ of
a geodesic ray whose point at infinity in $\PP_1(K)
\subset \partial_\infty \TT_\infty$ is equal to $\wt z$. Conversely,
any geodesic ray whose point at infinity lies in $\PP_1(K)
\subset \partial_\infty \TT_\infty$ contains a subray that maps
injectively by the canonical projection $\TT_\infty\ra \Ga_\infty \bs
\TT_\infty$.

Let us denote by $\widehat{\Ga_\infty\bs\TT_\infty}= (\Ga_\infty\bs
\TT_\infty)\sqcup \E_\infty$ Freudenthal's compactification (see
\cite{Freudenthal31}) of $\Ga_\infty \bs\TT_\infty$ by its finite set
of ends $\E_\infty$. This set of ends is indeed finite, in bijection
with $\Ga_\infty\bs\PP_1(K)$ by the map which associates to $z\in
\Ga_\infty\bs\PP_1(K)$ the end towards which the cuspidal ray $h_z$
converges. See for instance \cite{Serre83} for a geometric
interpretation of $\E_\infty$ in terms of the curve $\bf C$.

Let $\widehat{X_\infty}= X_\infty \sqcup \E_\infty$ and let
$\widehat{p_\infty}: \widehat{X_\infty} \ra \widehat{\Ga_\infty \bs
  \TT_\infty}$ be the map equal to the identity map on $\E_\infty$ and
to the canonical projection
$$
p_\infty: X_\infty=\Ga_\infty\bs G_\infty \ra 
\Ga_\infty\bs V\TT_\infty=\Ga_\infty\bs G_\infty/\uG(O_\infty)
$$ 
on $X_\infty$ (see the diagram at the beginning of Section
\ref{sec:background}). Since $p_\infty$ is a proper map, this defined
a compactification of $X_\infty$, by endowing $\widehat{X_\infty}$
with the compact metrisable topology generated by the open subsets of
$U$ and the sets $\widehat{p_\infty}^{-1}(U)$ with $U$ an open
neighborhood of a point in $\E_\infty$. We will say that $\E_\infty$
is the {\it set of cusps} of $X_\infty$, and we will indicate towards
which cusp of $X_\infty$ the escape of mass occurs.

\medskip For every $x\in X_\infty$, define the {\it height} of $x$ in
$X_\infty$ by
\begin{equation}\label{eq:defiheight}
\heit_\infty(x)= 
d_{\Ga_\infty\bs\TT_\infty}(p_\infty(x),\Ga_\infty \,*_\infty)\;.
\end{equation}
For every cusp $z\in {\cal E}_\infty$ of $X_\infty$, define the {\it
  height of $x$ in $X_\infty$ relative to the cusp $z$} by
$\heit_{\infty, \,z} (x) = 0$ if $p_\infty(x)$ does not belong to $h_z
(]0, +\infty[)$, and 
$$
\heit_{\infty, \,z}(x)= 
d_{\Ga_\infty \bs \TT_\infty} (p_\infty(x), h_z(0))
$$ 
otherwise.

\blemm \label{lem:deformheight}
For all $g'\in G_\infty$ and $x\in X_\infty$, we have
$$
|\heit_\infty(x)-\heit_\infty(xg')|\leq d_{\TT_\infty}(*_\infty,g'\,*_\infty)\;,
$$
and $|\heit_{\infty,\,z}(x)-\heit_{\infty,\,z}(xg')|\leq d_{\TT_\infty}
(*_\infty,g'\,*_\infty)$ for every cusp $z\in\E_\infty$ of $X_\infty$.  
\elemm

\dem Let $g\in G_\infty$ be such that $x=\Ga_\infty g$. We have
$p_\infty(x)=\Ga_\infty \,g\,*_\infty$ and $p_\infty(xg')=\Ga_\infty
\,g\,g'\,*_\infty$. By the triangle inequality and since the
projection map $\TT_\infty\ra \Ga_\infty\bs \TT_\infty$ does not
increase the distances, we have
\begin{align*}
|\heit_\infty(x)-\heit_\infty(xg')| &\leq 
d_{\Ga_\infty\bs\TT_\infty}(\Ga_\infty\,g\,*_\infty,\Ga_\infty \,g\,g'\,*_\infty)
\\ &\leq d_{\TT_\infty}(g\,*_\infty,g\,g'\,*_\infty)
=d_{\TT_\infty}(*_\infty,g'\,*_\infty)\;.
\end{align*}
The second assertion follows if $p_\infty(x)$ and $p_\infty(xg')$
simultaneously belong or do not belong to (the image of) $h_z$. If
for instance $p_\infty(x)$ belongs to $h_z$ and $p_\infty(xg')$ does
not belong to $h_z$, then 
$$
d_{\Ga_\infty \bs \TT_\infty} (p_\infty(x), h_e(0)) \leq 
d_{\Ga_\infty \bs \TT_\infty} (p_\infty(x), p_\infty(xg'))\,
$$ 
and the result holds as above. \cqfd

\bigskip \noindent{\bf Example: } Assume that $\bf C$ is the
projective line over $\FF_q$ and that $\infty$ is its usual point at
infinity. Then the (image of the) geodesic ray in $\TT_\infty$
starting from $*_\infty$ with point at infinity $\infty\in
\PP_1(K_\infty)$, which is
$$
n\in\NN\mapsto [O_\infty \times \pi_\infty^{n}\;O_\infty]=
a_\infty^{n}\;*_\infty \in V\TT_\infty\;,
$$ 
is a (weak) fundamental domain for the action of $\Ga_\infty$ on
$V\TT_\infty$: it injects onto $\Ga_\infty\bs V\TT_\infty$ by the
canonical map $\TT_\infty\ra \Ga_\infty\bs \TT_\infty$.

Hence $G_\infty=\coprod_{n\in\NN}\Ga_\infty \;a_\infty^{n}\; \uG(O_\infty)$.
For every $g\in G_\infty$, the height of $x=\Ga_\infty g$ is the unique
$n\in\NN$ such that $g\in \Ga_\infty \;a_\infty^{n}\; \uG(O_\infty)$.
Note that if one writes $g$ in the Cartan decomposition of $G_\infty$ as
$g\in \uG(O_\infty)\;a_\infty^{m}\; \uG(O_\infty)$ for some $m\in
\NN$, then $m=d_{\TT_\infty}(*_\infty,g\,*_\infty)\geq \heit_\infty(x)$, with
usually strict inequality.

The quotient graph of finite groups $\Ga_\infty\bs\!\bs \TT_\infty$, whose
underlying graph is the geodesic ray $\Ga_\infty\bs \TT_\infty$, is called
the {\it modular ray}. With $F_0=\uG(k_\infty)$, $F'_0=F_0\cap F_1$ and
$F_n=\big\{\begin{bmatrix} a & b \\ 0 &  d\end{bmatrix}\in
\Ga_\infty\;:\;v_\infty(b)\geq -n\big\}$, the modular ray $\Ga_\infty\bs\!\bs
\TT_\infty$ (which has only one end) is given by the following figure.

\begin{center}
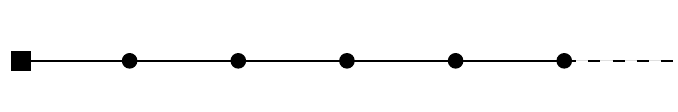
\end{center}

\bigskip
\noindent{\bf The full-down property in the modular graph (see for
  instance \cite{Serre83,Lubotzky91}). }
If $\rho$ is a geodesic ray in $\TT_\infty$ whose image is a cuspidal
ray in $\Ga_\infty\bs \TT_\infty$, the stabilizers of the vertices of
$\rho$ different from the origin of $\rho$ are strictly increasing
along the ray. Hence the image in $\Ga_\infty\bs \TT_\infty$ of a
geodesic ray in $\TT_\infty$ satisfies the following {\it full-down
  property}: if it starts to go down along the image of a cuspidal ray
$h_z$ for some $z\in\E_\infty$, then it needs to go all the way down
to $h_z(0)$.

As explained in \cite{Serre83,Paulin02}, this full-down property has
the following consequence: the image by the canonical map
$\TT_\infty\ra \Ga_\infty\bs \TT_\infty$ of a geodesic ray $\rho$ in
$\TT_\infty$ starting from $*_\infty$ either is an infinite sequence
$a_0b_0a_1b_1a_2b_2 \dots$ of concatenations of paths $a_i$ (possibly
reduced to points) in the finite graph $\Ga_\infty \bs
\TT_\infty-\bigcup_{z\in\E_\infty} h_z(]0,+\infty[)$ and back and
forth paths $b_i$ (of even lengths at least $2$) from the origin
$h_{z_i}(0)$ of the cuspidal ray $h_{z_i}$ to itself inside this ray,
if $\rho$ ends in an irrational point at infinity (that is, in
$\PP_1(K_\infty)-\PP_1(K)$), or starts by such a finite sequence and
then follows some cuspidal ray to infinity, otherwise.

\begin{center}
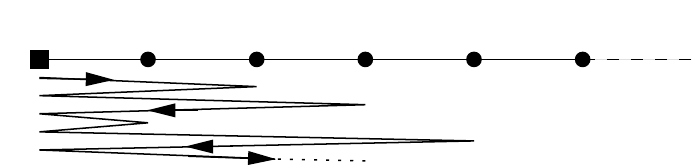
\end{center}

\subsection{$A_\infty$-periodic orbits in $X_\infty$}
\label{subsec:periodicorbit}

Let $(K,\nu,\uG,\uA)$ and the associated notation be as in Subsection
\ref{subsec:algebraic groups} before Remark \ref{rem:extension}.

Let us give a description of the compact orbits for the action by
translations on the right of the subgroup $A_\infty$ on
$X_\infty=\Ga_\infty\bs G_\infty$.

\bprop \label{prop:periodiccartanorbit} For every $g \in G_\infty$, the
following assertions are equivalent, where $x=\Ga_\infty g \in X_\infty$:
\begin{enumerate}
\item[(1)] there exists a unique $A_\infty$-invariant probability
  measure on the orbit $x A_\infty$;
\item[(2)] the subgroup $A_\infty\cap g^{-1} \Ga_\infty g$ is a
  (uniform) lattice in $A_\infty$;
\item[(3)] the orbit of $\pi'_\infty(x)$ under the geodesic flow
  $(\phi_n)_{n\in\ZZ}$ on $\Ga_\infty\bs\G\TT_\infty$ is periodic;
\item[(4)] there exists $\ga_0\in \Ga_\infty$ and $t_0\in
  K_\infty^\times$ with $v_\infty(t_0)$ positive and minimal such that
  $\ga_0 \,g= g \,\alpha_\infty(t_0)$.
\end{enumerate}  
\eprop

If one of these conditions is satisfied, we say that $x$ is {\it
  $A_\infty$-periodic}, and the unique $A_\infty$-invariant
probability measure on $x A_\infty$ is denoted by $\mu_x$. 

The elements $\ga_0\in\Ga_\infty$ and $t_0\in K_\infty^\times$ are
said to be {\it associated with} $g$. Note that they depend on the
choice of the representative $g$ of $x$: if $\ga_0$ is associated with
$g$, then $\ga^{-1}\ga_0\ga$ is associated with $\ga g$ for every
$\ga\in\Ga_\infty$. Furthermore, $\ga_0$ is primitive (not a proper
power of an element of $\Ga_\infty$) and loxodromic on
$\TT_\infty$. The period of $\pi'_\infty(x)$ under the geodesic flow
$(\phi_n)_{n\in\ZZ}$ is the translation length of $\ga_0$ on
$\TT_\infty$, which is equal to $v_\infty(t_0)$, and depends only on
$x$.

\medskip
\dem
The equivalence of (1) and (2) is well-known.

The equivalence of (2) and (3) follows from the equivariance of the
canonical projection $\pi'_\infty: X_\infty= \Ga_\infty \bs
G_\infty\ra \Ga_\infty\bs \G\TT_\infty=\Ga_\infty\bs G_\infty /
\uA(O_\infty)$ with respect to the morphism $v_\infty:A_\infty\ra\ZZ$ (see
Equation \eqref{eq:equivarvomega}).

The image $\Ga_\infty\,\ell= \Ga_\infty\, g\,\uA(O_\infty)$ in
$\Ga_\infty\bs\G\TT_\infty$ of the geodesic line
$\ell=g\uA(O_\infty)\in \G\TT_\infty= G_\infty/ \uA(O_\infty)$ is
periodic under the geodesic flow if and only if there exist $n>0$ and
$\ga_0\in\Ga_\infty$ such that $\ga_0\ell= \phi^n(\ell)=\ell
\,\alpha_\infty(\pi_\infty^n)$, hence, since
$\alpha_\infty:O_\infty^\times\ra\uA(O_\infty)$ is an isomorphism,
if and only if there exist $n>0$, $u_0\in O_\infty^\times$ and
$\ga_0 \in\Ga_\infty$ such that $\ga_0g=g \,
\alpha_\infty(\pi_\infty^n) \alpha_\infty(u_0)$. With
$t_0=\pi_\infty^nu_0$ so that $v_\infty(t_0)=n>0$, this proves the
equivalence of (3) and (4). \cqfd

\medskip 
Let us now prove the additional properties of $(\ga_0,t_0)$ and
discuss its uniqueness.  Assume that $n$ in the above proof is
minimal. Then $\ga_0$ is primitive and loxodromic, with translation
axis the image of $\ell$, translation length $n$, which is the period
of $\Ga_\infty\, \ell$ under the geodesic flow. Assume that
$(\ga'_0,t'_0) \in \Ga_\infty\times K_\infty^\times$ satisfies
$\ga'_0 \,g= g \,\alpha_\infty(t'_0)$ with $n'=v_\infty(t'_0)$
positive and minimal.  Then $n'=n$ and $\ga'_0\ell= \phi_{n}(\ell)$.
Hence $\ga_0^{-1}\ga'_0$ belongs to the pointwise stabilizer in
$\Ga_\infty$ of the image of $\ell$, which is the finite group
$g\uA(O_\infty)g^{-1}\cap \Ga_\infty$.  Therefore there exists
$u'_0\in {\alpha_\infty}^{-1} (g^{-1}\Ga_\infty g\cap\uA(O_\infty))
\subset O_\infty^\times$ such that $\ga'_0= \ga_0g
\alpha_\infty(u'_0) g^{-1}$ and $t'_0=t_0u'_0$.

\subsection{Hecke trees}
\label{subsec:Hecketrees}

Let $(K,\nu,\uG,\uA)$ and the associated notation be as in
Subsection \ref{subsec:algebraic groups} before Remark
\ref{rem:extension}.

The set $X_\infty$ of homothety classes of $R_\infty$-lattices in
$K_\infty\times K_\infty$ is the set of vertices of a graph, whose
non-oriented edges are the pairs $\{x,x'\}$ of vertices such that
there exists representatives $\Lambda$ of $x$ and $\Lambda'$ of $x'$
such that $\Lambda \subset \Lambda'$ and $\Lambda'/\Lambda$ is
isomorphic to $R_\infty/\pi_\nu R_\infty$. The action of $G_\infty$
on $X_\infty$ extends to an (isometric) action by graph automorphisms
on this graph.

For every $x\in X_\infty$, the connected component of $x$ in this
graph is a $(|k_\nu|+1)$-regular tree, called the {\it ($\nu$)-Hecke
  tree} of $x$, and denoted by $T_\nu(x)$. We have
$T_\nu(x)g=T_\nu(xg)$ for all $x\in X_\infty$ and $g\in G_\infty$. A
{\it ($\nu$-)Hecke ray} from $x$ is a geodesic ray in the Hecke tree
$T_\nu(x)$ starting from $x$.

\medskip The following description of the $\nu$-Hecke trees in
$X_\infty$ is well known, and is given, besides in order to fix the
notation, only for the sake of completeness.

\blemm \label{lem:relatheckbruh} Let $g\in G_\infty$ and $x=\Ga_\infty
\,g$ its image in $X_\infty$. The map from $G_\nu$ to $X_\infty$
defined by $g'\mapsto \pi_\infty(\Ga_S(g,g'))$ induces an isometric
map $\hec_g$ from the vertex set $V\TT_\nu= G_v/ \uG(O_\nu)$ of the
Bruhat-Tits tree $\TT_\nu$ onto the vertex set $VT_\nu (x)$ of the
Hecke tree $T_\nu(x)$, sending $*_\nu$ to $x$. For every
$\ga_0\in\Ga_\infty$, the map $\hec_g$ conjugates the action of
$\ga_0$ on $\TT_\nu$ to the right action of $g^{-1}\ga_0
g\in\Ga_\infty$ on $VT_\nu (x)$: for every $y\in V\TT_\nu$, we have
\begin{equation}\label{eq:commuthec}
  \hec_g(\ga_0 \,y)= \hec_g(y)\,g^{-1}\ga_0 g\;. 
\end{equation}
For all $h \in G_\infty$ such that $\Ga_\infty\, h=x$, we have
$\hec_g=\hec_h$ if and only if $g=h$; furthermore, the following
diagram commutes:
\begin{equation}\label{eq:diagheccommut}
\begin{array}[b]{ccc} 
V\TT_\nu & \stackrel{gh^{-1}}{\longrightarrow} & V\TT_\nu \\
\;_{\hec_h}\!\!\searrow & & \swarrow_{\hec_g}\\
 & VT_\nu(x)&
\end{array}.
\end{equation}
\elemm

Note that $\hec_g$ depends on $g$ and not only on $x$. We will denote
again by $\hec_g$ the (continuous) extension $\partial_\infty\TT_\nu
\ra \partial_\infty T_\nu(x)$ of $\hec_g$ to the
boundaries at infinity of the Bruhat-Tits and Hecke trees.

\medskip 
\dem 
Since the action by translations on the right of $\uG(O_\nu)$ on $X_S$
preserves the fibers of the bundle map $\pi_\infty:X_S\ra X_\infty$,
the map $g'\mapsto \pi_\infty(\Ga_S(g,g'))$ does induce a map
$\hec_g:V\TT_\nu= G_v/ \uG(O_\nu)\ra X_\infty$.

By definition of the Hecke tree $T_\nu(x)$ of $x=\Ga_\infty \,g=
g^{-1} [R_\infty\times R_\infty]$, its vertices are the points
$g^{-1}\ga \,[R_\infty\times \pi_\nu^n R_\infty]$ where
$\ga\in\Ga_\infty$ and $n\in\NN$.  By Equation
\eqref{eq:decomradialGnu}, any element in $G_\nu$ may be written
$\ga\,a_\nu^{n}\,g'$ for some $\ga\in\Ga_\infty$, $n\in\NN$ and
$g'\in \uG(O_\nu)$. Hence, the elements in $\hec_g(V\TT_\nu)$ are the
points $\pi_\infty(\Ga_S(g,\ga\,a_\nu^{n}\,g'))=\Ga_\infty
\,a_\nu^{-n}\, \ga^{-1} g$ where $g'\in \uG(O_\nu)$,
$\ga\in\Ga_\infty$ and $n\in\NN$.  Therefore $\hec_g(V\TT_\nu)
=VT_\nu(x)$.

If $y,y'\in V\TT_\nu$ are joined by an edge in $\TT_\nu$, then again
by density of $\Ga_\infty$ in $\uG(O_\nu)$, there exists an element in
$\Ga_\infty$ mapping the edge between $y$ and $y'$ into the geodesic
ray with vertices $(a_\nu^{n}\,*_ \nu)_{n\in\NN}$. Up to exchanging
$y$ and $y'$, there hence exists $n\in\NN$ and $\ga\in \Ga_\infty$
such that $\ga^{-1}y=\,a_\nu^{n}\,*_ \nu$ and $\ga^{-1}y'
=\,a_\nu^{n+1}\,*_ \nu$. In particular, $\hec_g(y)= \Ga_\infty
\,a_\nu^{-n}\,\ga^{-1}g$ is joined by an edge to $\hec_g(y')
=\Ga_\infty \,a_\nu^{-n-1} \,\ga^{-1}g$ in the Hecke tree
$T_\nu(x)$.  Hence $\hec_g$ induces a surjective graph morphism
between the trees $\TT_\nu$ and $T_\nu(x)$.  Since both trees are
regular of degree $|k_\nu|+1$, the map $\hec_g$ is an isomorphism of
trees.

Equation \eqref{eq:commuthec} follows by writing $y\in V\TT_\nu=
G_\nu/\uG(O_\nu)$ as $y=g'\uG(O_\nu)$ for some $g'\in \Ga_S$
(see the line following Equation \eqref{eq:decomradialGnu}), and by
using the following equalities:
\begin{align*}
\pi_\infty(\Ga_S(g,g'))g^{-1}\ga_0^{-1} g &=
\pi_\infty(\Ga_S({g'}^{-1}g,e))g^{-1}\ga_0^{-1} g=
\Ga_\infty({g'}^{-1}g)\;g^{-1}\ga_0^{-1} g\\ &=
\pi_\infty(\Ga_S(g,\ga_0g'))\;.
\end{align*}

Let $h$ be another element in $G_\infty$ such that $\Ga_\infty\,
h=x$. Since $G_\nu= \Ga_S\; \uG(O_\nu)$ and by the definition of
$\pi_\infty$, we have $\hec_g=\hec_h$ if and only if $\Ga_\infty
\ga^{-1}g= \Ga_\infty \ga^{-1} h$ for every $\ga\in\Ga_S$, that is
$\ga^{-1} (gh^{-1}) \ga\in\Ga_\infty$ for every $\ga\in\Ga_S$. Writing
$gh^{-1} = \begin{bmatrix}a & b \\ c & d \end{bmatrix}$ and using
$\ga=e$, we may take $a,b,c,d\in R_\infty$. Since the order of
vanishing at a point of ${\bf C}-\{\infty\}$ of an element of
$R_\infty$ is nonnegative and $v_\nu(\pi_\nu)=1$, we have
$v_\infty(\pi_\nu)\neq 0$ by the product formula. Taking
$\ga=\begin{bmatrix} \pi_\nu^n & 0 \\ 0 & 1
\end{bmatrix}$ gives $\pi_\nu^n c, \pi_\nu^{-n} b\in R_\infty$ for every 
$n\in\ZZ$, that is $c=b=0$. Taking $\ga=\begin{bmatrix}1 & \pi_\nu^n 
\\  0 & 1 \end{bmatrix}$ gives $\pi_\nu^n (a-d)\in R_\infty$ for every
$n\in\ZZ$, that is $a=d$. Hence $\hec_g=\hec_h$ if and only if
$gh^{-1}$ is the identity element in $\Ga_\infty=\uG(R_\infty)$.

The other claims are left to the reader. 
\cqfd

\section{Dynamics of the modular group at the infinite 
place on the Bruhat-Tits tree at a finite place}
\label{sec:dynaminftynu} 

Let $(K,\nu,\uG,\uA)$ and the associated notation be as in
Subsection \ref{subsec:algebraic groups} before Remark
\ref{rem:extension}.

In this Section, we study the dynamics of $\Ga_\infty$ on the
Bruhat-Tits tree $\TT_\nu$ of $(\uG,K_\nu)$.  Since $R_\infty\subset
O_\nu$, the lattice $\Ga_\infty=\uG(R_\infty)$ is contained in the
stabilizer $\uG(O_\nu)$ in $G_\nu$ of the base point $*_\nu$ in
$\TT_\nu$. Hence $\Ga_\infty$ does act on $\TT_\nu$, and for every
$n\in\NN$, every $\ga_0\in \Ga_\infty$ preserves the sphere
$$
S_\nu(n)=S_{\TT_\nu}(*_\nu,n)
$$ 
of center $*_\nu$ and radius $n$ in $\TT_\nu$. Since $S_\nu(n)$ is
finite, every orbit in $S_\nu(n)$ of the cyclic group ${\ga_0}^\ZZ$
generated by $\ga_0$ is periodic. The following linear growth property
of these periodic orbits is a remarkable feature of the positive
characteristic.

\btheo \label{theo:lineargrowpspherorbit} Let $\ga_0$ be an
element in $\Ga_\infty$ which is loxodromic on $\TT_\infty$. Let $\wt
K_\nu=\wt K_\nu(\ga_0)$ be the splitting field of $\ga_0$ over
$K_\nu$, with local ring $\wt O_\nu$, uniformizer $\wt \pi_\nu$
and residual field $\wt k_\nu$.   Let $e_\nu=e_\nu(\ga_0)$ be the
ramification index $e(\wt K_\nu,K_\nu)$ of $\wt K_\nu$ over $K_\nu$.
Let $d_\nu= d_\nu(\ga_0)$ be the smallest positive integer such that
the image of ${\ga_0}^{d_\nu}$ in $\uG(\wt k_\nu)$ (by reduction
modulo $\wt\pi_\nu \wt O_\nu$) is the identity. Let
$r_\nu=r_\nu(\ga_0)$ be the biggest positive integer such that the
image of ${\ga_0}^{d_\nu}$ in $\uG(\wt O_\nu/ {\wt
  \pi_\nu}^{r_\nu+1}\wt O_\nu)$ is not the identity. Then there
exists a constant $\kappa_\nu=\kappa_\nu(\ga_0) \in\NN$ such that for
every big enough $n\in \NN$, the maximal cardinality $m_n=m_n(\ga_0)$
of an orbit of ${\ga_0}^\ZZ$ in $S_\nu(n)$ satisfies
$$
m_n\leq e_\nu\;d_\nu\;
p^{\lceil\log_p\frac{n+\kappa_\nu}{r_\nu}\rceil}\;.
$$ 
\etheo

This result implies that the sequence $(m_n)_{n\in\NN}$ has linear
growth: for every $n\in\NN$ big enough, we have
\begin{equation}\label{eq:lineargrowth}
m_n\leq \frac{e_\nu\;d_\nu\;p}{r_\nu}\;(n+\kappa_\nu)\;,
\end{equation}
and that if $\ga_0$ is diagonalisable over $K_\nu$, then for every
$k\in\NN$ big enough
$$
m_{r_\nu p^k+\kappa_\nu} \leq d_\nu \,p^k\;.
$$

\medskip
\dem We start the proof by the following lemma on the growth of the
valuations of the powers of the elements of $O_\nu$ with their
constant terms removed, which concentrates the positive characteristic
feature.

\blemm \label{lem:onedimcase} 
Let $a\in k_\nu^{\;\times}$, $\lambda\in a+\pi_\nu O_\nu$ and $n\in
\NN$.  Define $m_n(\lambda) =\min\{k\in\NN-\{0\}\;:\; \lambda^k \in
a^k+\pi_\nu^n O_\nu\}$ and $r_\lambda=v_\nu(\lambda -a)>0$. Then
for every $n> r_\lambda$,
$$
m_n(\lambda)=p^{\lceil\log_p\frac{n}{r_\lambda}\rceil}\;.
$$
\elemm

In particular, $m_n(\lambda)< \frac{p}{r_\lambda}\;n$ for every
$n>r_\lambda$ and $m_{r_\lambda p^k}(\lambda)= p^k$ for every
$k\in\NN-\{0\}$.

\medskip
\dem Up to replacing $\lambda$ by $\frac{\lambda}{a}$, we may assume
that $a=1$. To simplify the notation, let $r=r_\lambda$.  For every
$k\in\NN-\{0\}$, consider the expansion of $k$ in base $p$ given by
$k=\sum_{i=0}^s a_i p^i$ where $s\in\NN$ and $a_i\in \{0,\dots,
p-1\}$.  Let $v_p(k)= \inf\{i\in\NN\;:\; \forall\;j<i,\;\; a_j=0\;\}$
be the $p$-adic valuation of $k$. Then, using the Frobenius
automorphism, and the fact that $a_i$ is invertible in the
characteristic subfield $\FF_p$, hence in $ O_\nu$, if and only if
$a_i$ is non zero, we have
$$
(1+{\pi_\nu}^r O_\nu^{\;\times})^k\subset
\prod_{i=0}^s (1+{\pi_\nu}^{rp^i} O_\nu^{\;\times})^{a_i}
\subset\prod_{0\leq i\leq s,\;a_i\neq 0} 
(1+{\pi_\nu}^{rp^i} O_\nu^{\;\times})
\subset 1+{\pi_\nu}^{rp^{v_p(k)}} O_\nu^{\;\times}\;.
$$
Hence for every $n\in\NN$, we have $\lambda^k\in 1+\pi_\nu^n O_\nu$ if
and only if $r p^{v_p(k)}\geq n$. Therefore, for every $n> r$, if
$rp^{m-1}< n \leq rp^{m}$ (that is, if $m= \lceil\log_p \frac{n}{r}
\rceil$), we have the equalities $m_n(\lambda)= \min\{k\in\NN-\{0\}
\;:\; v_p(k)=m\}= p^m$.  The result follows.  \cqfd

\medskip Now, let $\ga_0\in\Ga_\infty$ be loxodromic on
$\TT_\infty$. Note that the constant $d_\nu$ is well defined since
$R_\infty\subset O_\nu\subset \wt O_\nu$.  As we have seen in Section
\ref{subsec:BruhatTitstrees}, there exist $\lambda_\pm$ in a finite
extension of $K$ such that the element
$\ga_0$ is conjugated to $\begin{bmatrix} \lambda_+ & 0 \\
  0 & \lambda_-\end{bmatrix}$ and $\wt{K_\nu}=
K_\nu(\frac{\lambda_+}{\lambda_-})$. Note that $\lambda_-$ and
$\lambda_+$ are distinct since $\ga_0$ is not the identity element.

Let $\wt \TT_\nu$ be the Bruhat-Tits tree of $(\uG, \wt K_\nu)$, and
$\wt *_\nu=e\,\uG(\wt O_\nu)$ its standard base point in $V\wt
\TT_\nu=\uG(\wt K_\nu)/\uG(\wt O_\nu)$.  The value group of (the
unique extension of) the valuation $v_\nu$ on ${\wt K_\nu}^\times$
contains the value group $\ZZ$ of the valuation $v_\nu$ on
$K_\nu^\times$ with index $e_\nu$. By the correspondance of the action
on the right of $\uA(\wt K_\omega)$ on $\uG(\wt K_\omega)/\uA(\wt
O_\omega)$ and the action of the geodesic flow on the geodesic lines
in $\wt \TT_\nu$, the sphere $S_\nu(n)$ of center $*_\nu$ and radius
$n$ in $\TT_\nu$ is naturally contained in the sphere $S_{\wt \TT_\nu}
(\wt*_\nu, e_\nu\,n)$ of center $\wt *_\nu$ and radius $e_\nu\,n$ in
$\wt \TT_\nu$, for every $n\in\NN$.  Therefore, up to replacing
$K_\nu$ by $\wt K_\nu$, we may assume that $\ga_0$ is diagonalisable
over $K_\nu$, and we prove that the cardinality of every orbit of
${\ga_0}^\ZZ$ in $S_\nu(n)$ is at most $d_\nu\; p^{\lceil\log_p
  \frac{n+\kappa_\nu}{r_\nu}\rceil}$ for every $n\in\NN$, for some
$\kappa_\nu\in\NN$.

Note that the coefficients $\lambda_\pm$ have absolute value $1$ in
$K_\nu$.  Indeed,  $\begin{pmatrix} \lambda_+ & 0 \\ 0 & \lambda_- 
\end{pmatrix}$ may be choosen to be conjugated to a representative of
$\ga_0$ in $\GL_2(R_\infty)$. Hence $\lambda_\pm$ satisfy an equation
$P(\lambda_\pm)=0$ with $P$ a monic quadratic polynomial with
coefficients in $R_\infty\subset O_\nu$. Therefore
${|\lambda_\pm|_\nu}^2 \leq \max\{|\lambda_\pm|_\nu,1\}$, so that
$|\lambda_\pm|_\nu \leq 1$, and equality holds by replacing $\ga_0$ by
its inverse. Hence $\lambda_\pm\in a_\pm + \pi_\nu O_\nu$ with
$a_\pm\in k_\nu^\times$.

By the finiteness of $k_\nu^\times$, there exists a smallest $d_\nu
\in \NN-\{0\}$ such that ${a_-}^{d_\nu}= {a_+}^{d_\nu}$.  Note that
$d_\nu$ coincides with the notation introduced in the statement of
Theorem \ref{theo:lineargrowpspherorbit}. Let 
\begin{equation} \label{eq:defrnu}
r_\nu= v_\nu(\big(\frac{\lambda_+} {\lambda_-} \big)^{d_\nu}-1)\;.
\end{equation}
Since $\ga_0$ is loxodromic on $\TT_\infty$, no power of $\ga_0$ is
the identity, hence $r_\nu>0$. Note that $r_\nu$ coincides with the
notation introduced in the statement of Theorem
\ref{theo:lineargrowpspherorbit}. Up to replacing $\ga_0$ by
${\ga_0}^{d_\nu}$, to modify $\lambda_\pm$ by a common multiple by an
element of $k_\nu^\times$, and to proving that $m_n(\ga_0) \leq
p^{\lceil\log_p\frac{n+\kappa_\nu}{r_\nu}\rceil}$ for some
$\kappa_\nu\in\NN$ and for $n$ big enough, we may assume that the
constant terms in $k_\nu^\times$ of $\lambda_\pm$ are equal to $1$,
so that $d_\nu=1$.

Since $\ga_0$ is diagonalisable over $K_\nu$, there exists $h\in
G_\nu$ such that $\ga_0=h \begin{bmatrix} \lambda_- & 0 \\ 0 &
  \lambda_+ \end{bmatrix}h^{-1}$.  Since $\lambda_- \neq\lambda_+$,
the centralizer $Z_{G_\nu} (\ga_0)$ of $\ga_0$ in $G_\nu$ is the
abelian group $h\,A_\nu h^{-1}$. Note that $h$ is well defined modulo
multiplication on the right by an element of $A_\nu$.

Let $\ell_0:n\mapsto a_\nu^n\;*_\nu$ be the geodesic line in
$\TT_\nu$ from $0\in\partial_\infty\TT_\nu$ to $\infty \in
\partial_\infty\TT_\nu$, through $*_\nu$ at time $n=0$, which is
pointwise fixed by $\uA(O_\nu)$. The group $A_\nu$ preserves
$\ell_0(\ZZ)$ and acts transitively on it.  Note that the projective
action of $\uA(O_\nu)$ on $\PP^1(K_\nu)$ fixes $0$ and $\infty$, and
acts transitively on $\pi_\nu^{-k}O_\nu^\times\subset \PP^1(K_\nu)$
for every $k\in\ZZ$.

The geodesic line $\ell=h\,\ell_0$ is pointwise fixed by
$h\,\uA(O_\nu)\,h^{-1}$. Up to multiplying $h$ on the right by an
element of $A_\nu$, we may assume that the closest point to $*_\nu$ on
(the image of) $\ell$ is $h\,*_\nu=\ell(0)$. Let $s_\nu= s_\nu(\ga_0)
\in \NN$ be the distance between $*_\nu$ and $h \,*_\nu$ in $\TT_\nu$
(see the picture below).

\begin{center}
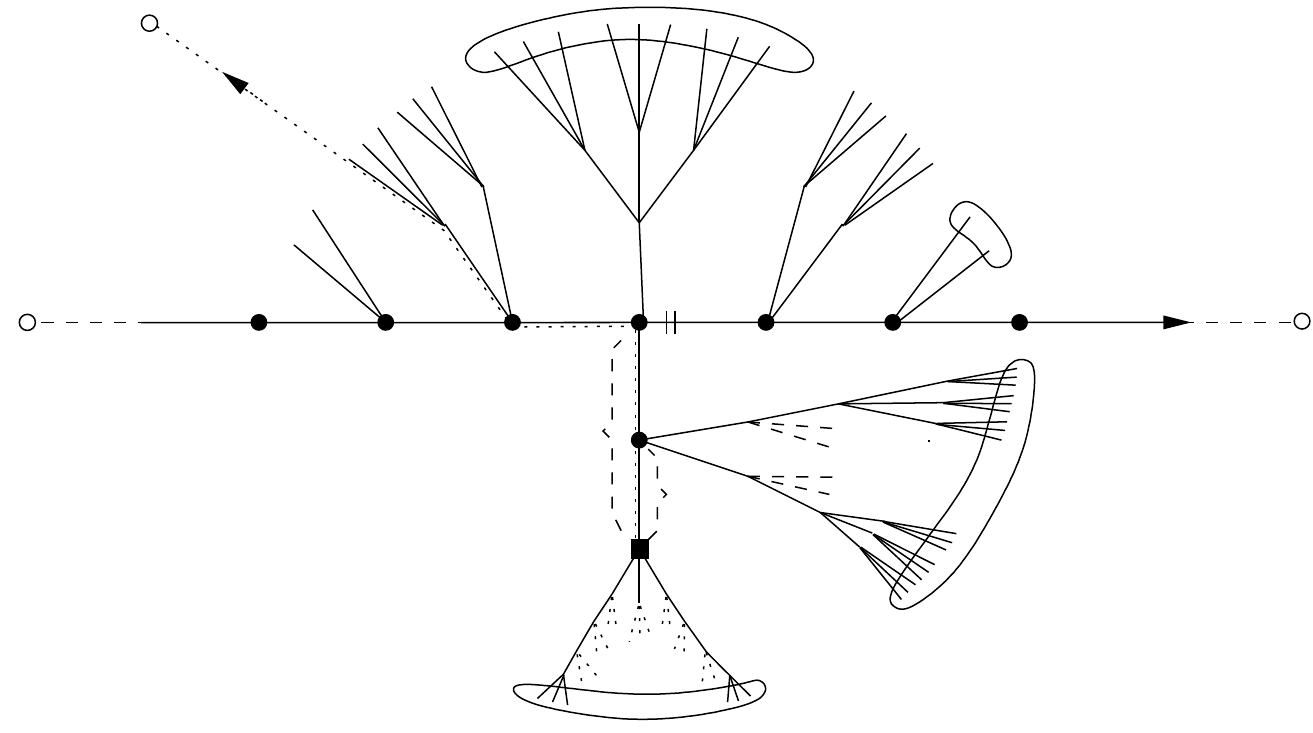
\end{center}

Let $\pr_{\ell}: V\TT_\nu\ra \ell(\ZZ)$ be the closest point map on
the geodesic line $\ell$.  For all $n,k\in\NN$, define (see the above
picture)
$$
E_{n,\,k}=\{x\in S_\nu(n)\;:\; \pr_{\ell}(x)=\ell(k)\}
$$
if $k\neq 0$, and 
$$
E_{n,\,0}=\big\{x\in S_\nu(n)\;:\; \pr_{\ell}(x)=\ell(0),\;\;\;
[h\,*_\nu,*_\nu]\cap[h\,*_\nu,x]= \{h\,*_\nu\}\,\big\}\;.
$$
For all $n,k'\in\NN$ with $0\leq k'<s_\nu$, let $E'_{n,\,k'}$ be the
set of $x\in S_\nu(n)$ such that the length of the common segment
$[*_\nu,h\,*_\nu]\cap[*_\nu,x]$ is equal to $k'$.
Then we have a partition
$$
S_\nu(n)=
\bigcup_{0\leq k'<s_\nu} E'_{n,\,k'}\cup\bigcup_{-n\leq k\leq n} E_{n,\,k}\;.
$$ 
Since $\ga_0$ fixes $*_\nu$, $h(0)$ and $h(\infty)$, it pointwise
fixes $\ell(\ZZ)\cup[*_\nu,h\,*_\nu]$. Hence the above partition of
$S_\nu(n)$ is invariant under $\ga_0$.

Note that $E_{n,\,k}$ is exactly the set of points at distance
$n-|k|-s_\nu$ from $ha_\nu^k*_\nu=\ell(k)$ on a geodesic ray from
$ha_\nu^k*_\nu$ to a point in $h(\pi_\nu^{-k} O_\nu^\times)\subset
\PP^1(K_\nu)$. Hence for any two points in $E_{n,\,k}$ (with $n,k$
fixed), there exists an element in the centralizer of $\ga_0$ mapping
one to the other. In particular, the cardinality $c_{n,\,k}=\card
(\ga_0^\ZZ y)$ is independent of $y\in E_{n,\,k}$.

Since $ha_\nu^{k-1}h^{-1}$ centralizes $\ga_0$ and $ha_\nu^{k-1}h^{-1}
E_{n-|k|+1,\,1}\subset E_{n,\,k}$, we have $c_{n,\,k}=
c_{n-|k|+1,\,1}$.  For every $n'\in\NN$, we have $c_{n',\,1} \leq
c_{n'+1,\,1}$, since the closest point map $E_{n'+1,\,1}\ra E_{n',\,1}$
is onto and equivariant under $\ga_0$.

Every point of $E'_{n,\,k'}$ is at distance $n+s_\nu-2k'$ from
$h\,*_\nu$. Hence $ha_\nu h^{-1}(E'_{n,\,k'})\subset E_{n+2s_\nu-2k'+1,\,1}$.
Therefore $c'_{n,\,k'}=\card(\ga_0^\ZZ y)$ is independent of $y\in
E'_{n,\,k'}$ and satisfies $c'_{n,\,k'}= c_{n+2s_\nu-2k'+1,\,1}$.

In particular, for every $n>s_\nu$, we have
$$
m_n(\ga_0)=\max\big\{\;\max_{0\leq k'<s_\nu} \;c'_{n,\,k'},
\;\max_{|k|\leq n}\; c_{n,\,k}\,\big\}=  c_{n+2s_\nu+1,\,1}\;.
$$

Note that $h(0)$ and $h(\infty)$ do not belong to $\PP_1(K)$, since
$\ga_0$, being loxodromic on $\TT_\infty$, fixes no point of
$\PP_1(K)$.  The positive subray of $\ell_0$ hence has no subray whose
image is entirely contained in the image of $\ell$. Therefore
$\ell_0(\NN)\cap \ell(\ZZ)$ is either empty or the set of vertices of
a compact interval $[\ell(0),\ell(k_0)]$ for some $k_0\in\ZZ$.

Assume first that $\ell_0(\NN)\cap \ell(\ZZ)$ is empty. Then the
segment $[*_\nu,h\,*_\nu]\cap[*_\nu,\infty[$ has length $k'_0\in
[0,s_\nu[\; \cap \,\NN$. Define $\kappa=2k'_0\in\NN$. Since the point
$a_\nu^{n'}*_\nu$ belongs to $E'_{n',\,k'_0}$ for every $n'\geq k'_0$,
the number $m_n(\ga_0)= c_{n+2s_\nu+1,\,1}=c'_{n+2k'_0,\,k'_0}$ is
the cardinality of the orbit under $\ga_0^\ZZ$ of
$a_\nu^{n+\kappa}*_\nu=[O_\nu \times {\pi_\nu}^{n+\kappa} O_\nu]$, if
$n$ is big enough.

Assume now that $\ell_0(\NN) \cap \ell(\ZZ)= [\ell(0),\ell(k_0)]\cap
V\TT_\nu$ for some $k_0\in\ZZ$ (see the above picture). Define $\kappa
= |k_0|+2s_\nu\in\NN$. Since the point $a_\nu^{n'}*_\nu$ belongs to
$E_{n',\,k_0}$ for every $n'\geq |k'_0|+s_\nu$, the number $m_n(\ga_0)
= c_{n+2s_\nu+1,\,1}=c_{n+|k_0|+2s_\nu,\,k_0}$ is the cardinality of
the orbit under $\ga_0^\ZZ$ of $a_\nu^{n+\kappa}*_\nu= [O_\nu \times
{\pi_\nu}^{n+\kappa} O_\nu]$, if $n$ is big enough.

For all $n\in\NN$, an element of $\GL_2(O_\nu)$ fixes $[O_\nu \times
\pi_\nu^{n} O_\nu]$ if and only if its $(2,1)$-coefficient vanishes
modulo $\pi_\nu^n$, that is, if it belongs to the Hecke congruence
subgroup of $\GL_2(O_\nu)$ modulo $\pi_\nu^n$.  Let
$\Ga_\infty(\pi_\nu^n)$ be the kernel of the morphism $\Ga_\infty\ra
\uG(R_\infty/\pi_\nu^nR_\infty)$ of reduction modulo ${\pi_\nu}^n$.
Thus for every $k\in\NN$, if $\ga_0^k$ belongs to
$\Ga_\infty({\pi_\nu}^{n+\kappa})$, then it fixes $\ell_0(n+\kappa)$.
Therefore, by the proof of Lemma \ref{lem:onedimcase} applied with
$\lambda=\frac{\lambda_+}{\lambda_-}$, since the constant $r_\lambda$
of Lemma \ref{lem:onedimcase} is equal to $r_\nu$ by Equation
\eqref{eq:defrnu}, we have, if $n$ is big enough,
\begin{align*}
m_n(\ga_0)&=\min\{k\in\NN-\{0\}\;:\; 
{\ga_0}^k\ell_0(n+\kappa)=\ell_0(n+\kappa)\}\\ &\leq
\min\{k\in\NN-\{0\}\;:\; \ga_0^k\in \Ga_\infty({\pi_\nu}^{n+\kappa })\}
\\ &\leq
\min\{k\in\NN-\{0\}\;:\; 
\;\; \big(\frac{\lambda_+}{\lambda_-}\big)^k\in 1+{\pi_\nu}^{n+\kappa} 
O_\nu\}\\ &= \min\{k\in\NN-\{0\}\;:\; 
\;\; v_p(k)\geq \log_p\frac{n+\kappa }{r_{\lambda}}\}\\ &=
\min\{k\in\NN-\{0\}\;:\; 
\;\; v_p(k)\geq \log_p\frac{n+\kappa }{r_\nu}\}=
p^{\lceil\log_p\frac{n+\kappa}{r_\nu}\rceil}\;.
\end{align*}
This concludes the proof of Theorem \ref{theo:lineargrowpspherorbit}. 
\cqfd

\section{Escape of mass along Hecke rays of
  $A_\infty$-periodic points}
\label{sec:ratheckeray}

Let $(K,\nu,\uG,\uA)$ and the associated notation be as in Subsection
\ref{subsec:algebraic groups} before Remark \ref{rem:extension}. We
fix from now on an $A_\infty$-periodic point $x_0$ in $X_\infty=
\Ga_\infty \bs G_\infty$, as well as a representative $g_0$ of $x_0$
in $\Ga_\infty$, so that $x_0=\Ga_\infty g_0$. In this section, we
prove our main results on the asymptotic behavior of the
$A_\infty$-invariant probability measures $\mu_x$ supported on the
$A_\infty$-orbits in $X_\infty$ of the vertices $x$ of the $\nu$-Hecke
tree $T_\nu(x_0)$ of $x_0$, as $x$ tends to infinity in this tree
along rays. We will recall below a proof that every vertex of
$T_\nu(x_0)$ is indeed $A_\infty$-periodic.

\medskip Let $\P(\widehat{X_\infty})$ be the space of probability
measures on the compactification $\widehat{X_\infty}= X_\infty\cup
\E_\infty$ by its finite set of cusps $\E_\infty=\Ga_\infty\bs
\PP_1(K)$ (see Subsection \ref{subsec:BruhatTitstrees}). Let $\xi
\in\partial_\infty T_\nu(x_0)$ be an end of the $\nu$-Hecke tree of
$x_0$. Let $\Theta_\xi$ be the subset of
$\P(\widehat{X_\infty})$ consisting of the weak-star accumulation
points of the sequence $(\mu_{x_n^\xi})_{n\in\NN}$ of $A_\infty$-invariant
probability measures on the vertices $(x^\xi_n)_{n\in\NN}$ along the
geodesic ray in $T_\nu(x_0)$ from $x_0$ to $\xi$.

For all $c>0$ and $z\in\E_\infty$, we say that
\begin{enumerate}
\item[$\bullet$] $\xi$ has {\it \cem} if there exists
  $\theta\in\Theta_\xi$ with  $\theta(\E_\infty)\geq c$.
\item[$\bullet$] $\xi$ has {\it \cem\ towards the cusp $z$} if there exists
  $\theta\in\Theta_\xi$ with $\theta(\{z\})\geq c$.
\item[$\bullet$] $\xi$ has {\it \ucem} if for every $\theta\in
  \Theta_\xi$ we have $\theta(\E_\infty)\geq c$.
\item[$\bullet$] $\xi$ has {\it \ucem\ towards the cusp $z$} if for
  every $\theta\in \Theta_\xi$ we have $\theta(\{z\})\geq c$.
\end{enumerate}

\subsection{Uniform escape of mass along rational 
Hecke rays}
\label{subsec:ratheckeray}

We start this subsection by defining the {\it rational Hecke rays} in
the $\nu$-Hecke tree $T_\nu(x_0)$ of $x_0$, and we will then prove
Theorem \ref{theo:lossmass}, a uniform escape of mass phenomenon for
the $A_\infty$-invariant probability measures $\mu_x$, as $x$ tends to
infinity along these rays.

\medskip
The group $\uG(K)$ acts transitively on $\PP_1(K)$, but its subgroups
$\Ga_\infty=\uG(R_\infty)$ and $\Ga_S=\uG(R_\infty[\pi_\nu^{-1}])$ do
not in general. The sets $\E_\infty=\Ga_\infty\bs \PP_1(K)$ (with
order at most the class number of $R_\infty$) and $\Ga_S\bs\PP_1(K)$
are finite and both canonical maps $\Ga_\infty\bs \PP_1(K)
\ra\Ga_S\bs\PP_1(K)\ra \uG(K)\bs\PP_1(K)$ may be non injective. Note
that for instance when $\bf C$ is the projective line over $\FF_q$ and
$\infty$ its usual point at infinity, then $R_\infty$ is principal,
and $\Ga_\infty$ does act transitively on $\PP_1(K)$.

Since $\Ga_\infty$ preserves $\PP_1(K)$ and by the commutativity of
the diagram \eqref{eq:diagheccommut}, the image $\hec_{g_0}(\PP_1(K))
\subset \partial_\infty T_\nu(x_0)$ by $\hec_{g_0}$ of the set
$\PP_1(K)$ of rational points of $\partial_\infty \TT_\nu=
\PP_1(K_\nu)$ does not depend on the choice of the representative
$g_0$ of $x_0$, nor does the image by $\hec_{g_0}$ of the orbit of
$\infty$ by any subgroup of $\uG(K)$ containing $\Ga_\infty$, as for
instance $\hec_{g_0}(\Ga_S\infty)$.

A Hecke ray in $T_\nu(x_0)$, as well as its point at infinity, is said
to be {\it rational} if its point at infinity belongs to
$\hec_{g_0}(\PP_1(K))$, and {\it $S$-rational} if its point at
infinity belongs to $\hec_{g_0}(\Ga_S\infty)$. In particular when
$\Ga_\infty$ acts transitively on $\PP_1(K)$ (that is, when the graph
$\Ga_\infty\bs \TT_\infty$ has only one end, as for instance when
${\bf C}$ is the projective line over $\FF_q$ and $\infty$ its usual
point at infinity), these two notions coincides. But there are
examples of functions fields when not all rational ends of
$T_\nu(x_0)$ are $S$-rational (the two inclusions $\Ga_\infty
\infty\subset \Ga_S\infty\subset \PP_1(K)$ may be strict).

If $\xi$ is a rational end of $T_\nu(x_0)$, the {\it cusp of
  $X_\infty$ associated to $\xi$} is $z_\xi=\Ga_\infty\ga \infty\in
\E_\infty$, where $\ga\in\uG(K)$ is such that $\xi=\hec_{g_0}(\ga
\infty)$. Note that $z_\xi$ does not depend on the choices of $g_0$ or
$\ga$. If $\xi$ is $S$-rational, we say that $z_\xi$ is an {\it
  $S$-cusp} of $X_\infty$.

\btheo \label{theo:lossmass} There exists $c=c(x_0)>0$ such that every
rational end $\xi$ of the Hecke tree of $x_0$ has \ucem, and if
furthermore $\xi$ is $S$-rational, then $\xi$ has \ucem\ towards the
cusp of $X_\infty$ associated to $\xi$.  
\etheo

\dem 
We start the proof by giving some notation. Let us fix elements
$\ga_0\in\Ga_\infty$ and $t_0\in K_\infty^\times$ associated with
the choosen representative $g_0$ of $x_0$ (see Proposition
\ref{prop:periodiccartanorbit} and its following comment): we have
$$
\ga_0 \,g_0= g_0 \,\alpha_\infty(t_0)
$$ 
and $\ell_0=v_\infty(t_0)>0$ is the translation distance of $\ga_0$ on
$\TT_\infty$.

Since $\uG(K)$ acts transitively on $\PP_1(K)$ and
$\E_\infty=\Ga_\infty\bs \PP_1(K)$ is finite, there exists a finite
subset $F_1$ of $\uG(K)$ such that $\PP_1(K)=\Ga_\infty F_1 \infty$, and
we may assume that $\Ga_S \infty= \Ga_\infty (F_1\cap \Ga_S)\infty$.

Since $\uG(K)$ commensurates $\Ga_S$, there exists a finite subset
$F_2$ of $\uG(K)$ such that for all $\ga\in F_1$ and $n\in\NN$, there
exists $b_{\ga,n}$ in $F_2$ such that
\begin{equation}\label{eq:defibgan}
\ga \,a_\nu^n \,\ga^{-1} \in \Ga_S \;b_{\ga,n}\;.
\end{equation}
We assume that $1\in F_2$ and $b_{\ga,n}=1$ if $\ga\in\Ga_S$. 

For every $b\in \uG(K)$, let $\ov b\in\Ga_S$ be such that $b\in \ov
b\,\uG(O_\nu)$, which exists by Equation \eqref{eq:decomradialGnu}.
We assume that $\ov b=1$ if $b\in\Ga_S$.

\medskip 
Now that this notation has been given, we consider the rational ends
$\xi$ of the Hecke tree $T_\nu(x_0)$. Let $\ga'=\ga'_\xi\in
\Ga_\infty$ and $\ga=\ga_\xi\in F_1$ be such that $\xi=\hec_{g_0}(\ga'\ga
\infty)$. We assume that $\ga\in\Ga_S$ if $\xi$ is $S$-rational.

\begin{center}
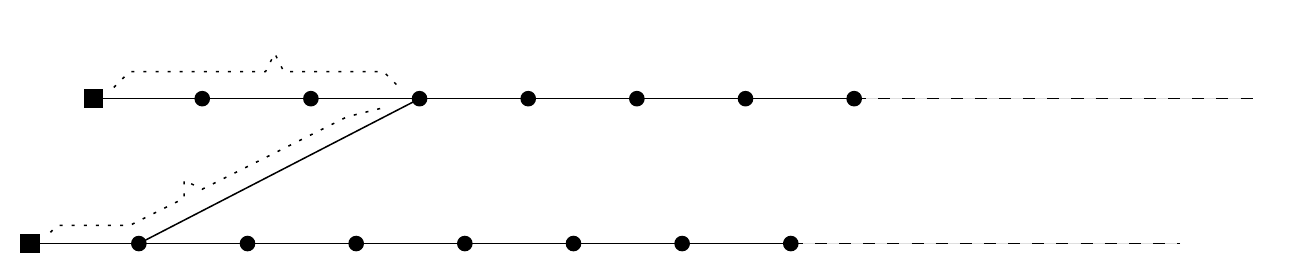
\end{center}

For every $n\in\NN$, let $y_n=a_\nu^n \,*_\nu$, so that for
every rational end $\xi$ of $T_\nu(x_0)$, the point at infinity of the
image by $\hec_{g_0}$ of the geodesic ray $n\mapsto \ga'\ga y_n$ is
$\xi$. Let $n'_\xi\in\NN$ be the distance from $*_\nu$ to this ray, and
let $n_\xi\in\NN$ be such that 
$$
[*_\nu,\ga'\ga \infty[\;\cap\,
[\ga'\ga*_\nu,\ga'\ga \infty[\;=[\ga'\ga y_{n_\xi},\ga'\ga \infty[\;.
$$
Let $r_\xi= n_\xi-n'_\xi\in \ZZ$.  Denote by $(x_n=
x_{n}^\xi)_{n\in\NN}$ the geodesic ray in the Hecke tree $T_\nu(x_0)$
from $x_0$ to $\xi$, for every $n\in\NN$ with $n\geq n_\xi$.
Using in the following sequence of equalities

$\bullet$~ the definition of $\hec_{g_0}$ for the third equality and 

$\bullet$~ the definition of $\pi_\infty$ (since $\ga
\,a_\nu^{n}\,\ga^{-1} b_{\ga,n}^{-1}\in\Ga_S$ by Equation
\eqref{eq:defibgan}) for the fifth one,

\noindent we have
\begin{align}
x_{n-r_\xi}& = \hec_{g_0}(\ga'\ga y_{n})=
\hec_{g_0}(\ga'\ga \,a_\nu^{n} \,*_\nu)=
\pi_\infty(\Ga_S({g_0},\ga'\ga \,a_\nu^{n}))\nonumber\\ & = \pi_\infty(
\Ga_S({g_0},\ga'\ga \,a_\nu^{n}\,\ga^{-1}\,b_{\ga,n}^{-1}\,b_{\ga,n}\ga))=
\Ga_\infty\,(\ga'\ga \,a_\nu^{n}\,\ga^{-1}b_{\ga,n}^{-1}
\,\overline{b_{\ga,n}\ga}\,)^{-1}{g_0}\nonumber\\ & =
\Ga_\infty\,((\,\overline{b_{\ga,n}\ga}\,)^{-1}b_{\ga,n}\ga)\,
a_\nu^{-n}\,(\ga{\ga'})^{-1}{g_0}\;.\label{eq:computxn}
\end{align}
Let $\E_\infty(\xi)$ be the subset of $\E_\infty$ consisting of the
elements
$$
z_n=\Ga_\infty(\, \overline{b_{\ga,n}\ga} \,)^{-1} b_{\ga,n}\ga \infty
$$ 
as $n$ varies. When $\xi$ is $S$-rational, we have $x_{n-r_\xi} =
\Ga_\infty\ga\, a_\nu^{-n}\,(\ga'\ga)^{-1}g{g_0}$, and $\E_\infty(\xi)$ is
the singleton of the cusp $z_\xi=\Ga_\infty\ga'\ga \infty=\Ga_\infty\ga
\infty$ associated to $\xi$.


Let $\pi'_\infty:X_\infty=\Ga_\infty\bs G_\infty\ra \Ga_\infty\bs
\G\TT_\infty=\Ga_\infty\bs G_\infty/\uA(O_\infty)$ be the canonical
projection, which is equivariant under $v_\infty: A_\infty \ra \ZZ$
(see Subsection \ref{subsec:periodicorbit}).  It extends to a
continuous map from $\widehat{X_\infty}= X_\infty \sqcup \E_\infty$ to
Freudenthal's compactification $\widehat{\Ga_\infty
  \bs\G\TT_\infty}=(\Ga_\infty \bs\G\TT_\infty) \sqcup \E_\infty$, by
the identity on $\E_\infty$.

By Equation \eqref{eq:commuthec}, since $\ga_0 \,{g_0}= {g_0}
\,\alpha_\infty(t_0)$ and by Equation \eqref{eq:equivarvomega}, for
every $k\in\NN$ and $y\in V\TT_\nu$, we have
\begin{align*}
\pi'_\infty(\hec_{g_0}({\ga_0}^ky)) &=
\pi'_\infty(\hec_{g_0}(y){g_0}^{-1}{\ga_0}^k{g_0})=
\pi'_\infty(\hec_{g_0}(y)\alpha_\infty({t_0}^k))\\ & =
\phi_{\ell_0 k}(\pi'_\infty(\hec_{g_0}(y)))\;.
\end{align*}
In particular, since the orbits of $\ga_0$ on $V\TT_\nu$ are finite,
every $x\in VT_\nu(x_0)$ is also $A_\infty$-periodic and the
$A_\infty$-invariant probability measure $\mu_{x}$ on the compact
orbit $x A_\infty$ is well defined. Furthermore, with the notation of
Theorem \ref{theo:lineargrowpspherorbit}, for every $n\geq n'_\xi$, the
orbit under the geodesic flow of $\pi'_\infty(x_n)= \pi'_\infty(
\hec_{g_0}(\ga'\ga y_{n+r_\xi}))$ is periodic, with period
$\lambda_n$ bounded as follows:
\begin{equation}\label{eq:majolengthorbit}
\lambda_n\leq \ell_0\,\min\big \{k\in\NN-\{0\}\;:\;
{\ga_0}^k\ga'\ga y_{n+r_\xi}=\ga'\ga y_{n+r_\xi}\big\} \leq
\ell_0\, m_{n}(\ga_0)\;.
\end{equation}

Let $d$ be the distance in the graph $\Ga_\infty\bs \TT_\infty$.
Recall that $p_\infty:X_\infty\ra \Ga_\infty\bs \TT_\infty$ is the map
$\Ga_\infty g\mapsto \Ga_\infty g *_\infty$ (see the diagram at the
beginning of Section \ref{sec:background}).  Using
\begin{enumerate}
\item[$\bullet$] Lemma \ref{lem:deformheight} with $\kappa= \kappa_\xi=
d_{\TT_\nu}(*_\nu,(\ga'\ga)^{-1} {g_0}\,*_\nu)$ for the first
inequality,
\item[$\bullet$] the definition of the height (see Equation
\eqref{eq:defiheight}) and Equation \eqref{eq:computxn} with the
notation $\beta_n=(\,\overline{b_{\ga,n+r_\xi}\ga}\,)^{-1}
b_{\ga,n+r_\xi}\ga$ for the second equality,
\end{enumerate}
\noindent 
we have, for all $n\in\NN$,
\begin{align}\label{eq:minoheightprelim}
  \heit_\infty(x_n) &\geq
  \heit_\infty(x_n(\ga^{-1}{\ga'}^{-1} {g_0})^{-1}) -\kappa \nonumber \\
  &= d(p_\infty(\Ga_\infty \beta_n \,a_\nu^{-n-r_\xi}),\Ga_\infty
  *_\infty) -\kappa \nonumber\\ &= d(\Ga_\infty
  \beta_n\,a_\nu^{-n-r_\xi}*_\infty,\Ga_\infty *_\infty) -\kappa\;.
\end{align}
Recall that $\beta_n$ belongs to $\uG(K)$, hence preserves the set of
geodesic rays in $\TT_\infty$ ending in
$\PP_1(K)\subset \partial_\infty\TT_\infty$, and takes finitely many
values as $n$ varies. Let 
$$
\kappa'=\kappa'_\xi=\max_{n\in\NN}\;
d(\Ga_\infty \beta_n*_\infty,\Ga_\infty*_\infty)+\kappa\;.
$$ 
Recall that any geodesic ray in $\TT_\infty$ ending in $\PP_1(K)$ has
a subray that isometrically injects into $\Ga_\infty\bs\TT_\infty$.
Hence using Equation \eqref{eq:minoheightprelim} and the triangle
inequality, there exist constants $n''_\xi\geq n'_\xi$ and
$\kappa''_\xi,\kappa'''_\xi\geq 0$ such that for every integer $n\geq n''_\xi$,
\begin{align}\label{eq:minoheight}
  \heit_\infty(x_n) & \geq d(\Ga_\infty
  \beta_n\,a_\nu^{-n-r_\xi}*_\infty,\Ga_\infty \beta_n*_\infty)
  -\kappa' \nonumber \\ & \geq
  d_{\TT_\infty}(a_\nu^{-n-r_\xi}*_\infty,*_\infty)-\kappa''_\xi \nonumber
  \\ &=(n-r_\xi)\,|v_\infty( \pi_\nu)| -\kappa''_\xi=n\,|v_\infty(
  \pi_\nu)|-\kappa'''_\xi\;.
\end{align}
This argument in fact proves that $\heit_{\infty,\,z_n}(x_n)\geq
n\,|v_\infty( \pi_\nu)|-\kappa'''_\xi$ for $n$ big enough, where $z_n$
is the cusp defined above.

For every $n\in\NN$, let $\mu'_n= (\pi'_\infty)_*\mu_{x_n}$, which is
the equiprobability on the finite orbit of $\pi'_\infty(x_n)$ under
the geodesic flow on $\Ga_\infty \bs \G\TT_\infty$. Recall that the
pushforwards of measures by proper continuous maps preserve the total
mass, and are weak-star continuous. The map $\pi'_\infty$ is a
fibration with compact fiber, hence a proper map. Therefore $\xi$ has
\ucem\ (repectively \ucem\ towards its associated cusp
$z_\xi=\Ga_\infty\ga'\ga \infty$) if and only if for every weak-star
accumulation point $\theta'$ of $(\mu'_{n})_{n\in\NN}$ in the space of
probability measures on $\widehat{\Ga_\infty \bs\G\TT_\infty}$, we
have $\theta'(\E_\infty)\geq c$ (respectively $\theta'(\{z_\xi\})\geq
c$).

Let $o:\Ga_\infty\bs\G\TT_\infty\ra \Ga_\infty\bs V\TT_\infty$ be the
origin map $\Ga_\infty\ell\mapsto \Ga_\infty \ell(0)$, which is a
proper map. For all $N\in\NN$, let 
$$
K_N=o^{-1}\big(\big\{x\in\Ga_\infty\bs
V\TT_\infty \;:\; x\in\bigcup_{z\in\E_\infty(\xi)}h_z([0,+\infty[),\;\;
d(x, \Ga_\infty *_\infty)\geq N\big\}\big)\;,
$$ 
which are open subsets of $\Ga_\infty\bs\G\TT_\infty$, which
accumulate as $N\ra+\infty$ exactly to $\E_\infty(\xi)\subset
\widehat{\Ga_\infty\bs \G\TT_\infty}$.  By the full-down property (see
Subsection \ref{subsec:BruhatTitstrees}), the orbit under the geodesic
flow of $\pi'_\infty(x_n)$ passes at a distance from $\Ga_\infty
*_\infty$ which is bounded by the diameter $N_0$ of the finite graph
$\Ga_\infty\bs\TT_\infty-\bigcup_{z\in\E_\infty}h_z(]0+\infty[)$. Recall
that this orbit is periodic, of period denoted by $\lambda_n$. Hence
if $N\geq N_0$ and if $\heit(x_n)\geq N$, the origins of $\phi_i
(\pi'_\infty(x_n))$ for $0\leq i\leq\lambda_n$ needs to range twice
over all points at distance between $N$ and $\heit(x_n)$ on a geodesic
ray in $\Ga_\infty\bs \TT_\infty$ between $\Ga_\infty *_\infty$ and
$o(\rho_\infty(x_n))$.  Hence if $n$ is big enough, by the comment
following Equation \eqref{eq:minoheight} and by Equation
\eqref{eq:majolengthorbit}, we have
\begin{equation}\label{eq:minomassatinfty}
  \mu'_n(K_N)\geq \frac{2(\heit_\infty(x_n)-N)}
{\lambda_n}\geq \frac{2n\,|v_\infty(\pi_\nu)|-2\kappa'''_\xi}
{\ell_0\, m_{n}(\ga_0)}\;.
\end{equation}
By the linear growth property of $(m_{n}(\ga_0))_{n\in\NN}$ (see
Equation \eqref{eq:lineargrowth} and the notation of Theo\-rem
\ref{theo:lineargrowpspherorbit}), the right hand side of Equation
\eqref{eq:minomassatinfty} as a limit as $n\ra+\infty$ at least
$$
c=\frac{2\,r_\nu(\ga_0)\,|v_\infty(\pi_\nu)|}
{\ell_0\;e_\nu(\ga_0)\;d_\nu(\ga_0)\;p}\;.
$$ 
Hence for every weak-star accumulation point $\theta'$ of
$(\mu'_{n})_{n\in\NN}$, we have $\theta'(\E_\infty(\xi))\geq
c$. This proves the result.  \cqfd

\medskip
\noindent {\bf Remark. } The aim of this remark is to give some
estimations on the constant $c$ appearing in this proof, and to give
examples of full escape of mass along rational Hecke rays.

As above, let $\ga_0\in\Ga_\infty$ be a (primitive loxodromic on
$\TT_\infty$) element associated to $x_0$, and let us fix $\wt\ga_0\in
\GL_2(R_\infty)$ whose image in $\Ga_\infty=\PGL_2(R_\infty)$ is
$\ga_0$.  Note that $\det\wt{\ga_0}\in R_\infty^\times=
k_\infty^\times$, hence $v_\infty(\det\wt{\ga_0})= 0$.  We may denote
by $\lambda_\pm$ the eigenvalues of $\wt{\ga_0}$ with
$v_\infty(\lambda_+)>0$, so that $v_\infty(\lambda_-)=
-v_\infty(\lambda_+) <v_\infty(\lambda_+)$ and, by Equation
\eqref{eq:translatlength},
$$
\ell_0= 2\, |v_\infty(\lambda_-)|=2\, |v_\infty(\tr(\wt{\ga _0} ))|\;.
$$ 
With the notation of Theorem \ref{theo:lineargrowpspherorbit}, let
us define
$$
\operatorname{LOM}(\ga_0)= \frac{r_\nu(\ga_0)\;|v_\infty(\pi_\nu)|}
{|v_\infty(\tr\wt{\ga_0})|\;e_\nu(\ga_0)\;d_\nu(\ga_0)}\;,
$$
so that we chose $c=\operatorname{LOM}(\ga_0)/p$ in the above proof.

Let us consider $n_k= r_\nu(\ga_0)\,p^k-\kappa_\nu(\ga_0)$ for
$k\in\NN$ big enough (again with the notation of Theo\-rem
\ref{theo:lineargrowpspherorbit}), so that $m_{n_k}(\ga_0) \leq
e_\nu(\ga_0)\,d_\nu(\ga_0)\,p^k$ by Theorem
\ref{theo:lineargrowpspherorbit}. Using this majoration on the
denominator in Equation \eqref{eq:minomassatinfty}, the above proof
gives moreover that every weak-star accumulation point $\theta'$ of
$(\mu'_{n_k})_{k\in\NN}$ satisfies $\theta'(\E_\infty(\xi))\geq
\operatorname{LOM}(\ga_0)$.  In particular, the sequence
$(\mu_{x_{n_k}})_{k\in\NN}$ weak-star converges to the $0$ measure on
$X_\infty$ if $\operatorname{LOM}(\ga_0) =1$. Let us give an example
of this when $\bf C$ is the projective line and $\infty$ its usual
point at infinity. Let $d=d_\nu(\ga_0)$, $e= e_\nu(\ga_0)$ and
$r=r_\nu(\ga_0)$, so that $\operatorname{LOM}(\ga_0)=
\frac{r\;|v_\infty(\pi_\nu)|} {e\;d\;|v_\infty(\tr\wt{\ga_0})|}$. Let
$\wt{k_\nu}$ be the residual field of the splitting field $\wt K_\nu$
of $\wt{\ga_0}$ over $K_\nu$.

\blemm Assume that the discriminant $\Delta=(\tr\wt{\ga_0})^2- 4\det
\wt{\ga_0}$ of $\wt{\ga_0}$ is irreducible over $\FF_q$, and let
$\pi_\nu=\Delta$. Then $\operatorname{LOM} (\ga_0) =1$.  
\elemm

This assumption is for instance satisfied if $-1$ is not a square
modulo $p$ (as for $p=3$), if $p=q$ and if $\wt\ga_0= \begin{pmatrix}
  Y & 1\\ 1 &0\end{pmatrix}$, since $\Delta=Y^2+4$. By the previous
arguments, for every rational end $\xi\in\Omega$, there exists an
element $\theta'\in\Theta_\xi$ which vanishes on $X_\infty$. This
proves Theorem \ref{theo:fullloss} in the introduction. The above
proof also gives a speed of escape of mass when $\operatorname{LOM}
(\ga_0) =1$: for every compact subset $C$ of $X_\infty$, we have
$\mu_{x_{n_k}}(C)=\bigO(\frac{1}{n_k})$ when $n_k=
r_\nu(\ga_0)\,p^k+\kappa_\nu(\ga_0)$.

\medskip \dem Since $\Delta$ is irreducible, we have $p\neq 2$. In
particular, the roots of $\wt\ga_0$ are $\lambda_\pm= \frac{1}{2}
\big(\tr\wt{\ga_0}\pm\sqrt{\pi_\nu}\,\big)$. We have $\wt K_\nu=
K_\nu(\sqrt{\pi_\nu}\,)$, and $\wt{k_\nu}=k_\nu$. In particular, the
ramification index of the splitting field of $\wt\ga_0$ over $K_\nu$
is $e=2$.  Since the constant terms in $\wt{k_\nu}$ (modulo
$\sqrt{\pi_\nu}$) of $\lambda_\pm$ are equal, we have $d=1$.  Since
$\deg(\det\wt{\ga_0}) =0$, we have
$$
|v_\infty(\pi_\nu)|=\deg\big((\tr\wt{\ga_0})^2- 4\det
\wt{\ga_0}\big)=2\deg(\tr\wt{\ga_0})=2\,|v_\infty(\tr\wt{\ga_0})|\;.
$$ 
Since $\wt\ga_0$ is not congruent to the identity modulo
$\sqrt{\pi_\nu}^{\,2}=\pi_\nu$, we have $r=1$. Hence
$\operatorname{LOM}(\ga_0) =1$.  \cqfd

\medskip 
Let us give one more estimation on the constant
$\operatorname{LOM}(\ga_0)$ when $p\neq 2$ and $v_\nu(\tr\wt{\ga_0})
>0$. We then have
$$
\lambda_\pm=\frac{1}{2}\big(\tr\wt{\ga_0}\pm
\sqrt{(\tr\wt{\ga_0})^2-4\det\wt{\ga_0}}\;\big).
$$
Since $\det\wt{\ga_0}\in k_\infty^\times \subset O_\nu^\times$, we
have $v_\nu(-4\det\wt{\ga_0}) = 0$, hence $e=1$ (and
$[\wt{k_\nu}:k_\nu]=1$ if $-\det\wt{\ga_0}$ is a square and $2$
otherwise). The constant terms $a_\pm=\pm\sqrt{-4\det\wt{\ga_0}}\in
\wt{k_\nu}^\times$ of $\lambda_\pm$ are opposite (and non zero), hence
$d=2$. By Equation \eqref{eq:defrnu}, we have $r=
v_\nu\big(\frac{\lambda_+}{\lambda_+}-1\big)=v_\nu(\tr\wt{\ga_0})$.
Furthermore
$$
|v_\infty(\tr\wt{\ga_0})| =\deg
(\tr\wt{\ga_0})\geq v_\nu(\tr\wt{\ga_0})\deg \pi_\nu=
r\, |v_\infty(\pi_\nu)|\;.
$$ 
Hence $\operatorname{LOM}(\ga_0)\leq \frac{1}{2}$, with equality if
and only if $\tr\wt{\ga_0}$ is a constant multiple of a power of
$\pi_\nu$, as for instance when $\pi_\nu=Y$ and
$\wt\ga_0=\begin{pmatrix} Y & 1\\ 1 &0\end{pmatrix}$. For these
elements where equality holds, at least half the mass escapes to
infinity along subsequences of every rational Hecke ray.

\subsection{Escape of mass along uncountably many
  Hecke rays}
\label{subsec:uncountableeom}

In the previous subsection, we proved escape of mass phenomena along
countably many Hecke rays, the rational ones. In this subsection, we
use the uniformity of the escape of mass in Theorem
\ref{theo:lossmass} in order to prove that an escape of mass (towards
prescribed cusps of $X_\infty$) actually occurs along uncountably many
Hecke rays. We first introduce some notation that we will use from now
on in this paper.

We denote by $\Om=\partial_\infty T_\nu(x_0)$ the boundary at infinity
of the Hecke tree $T_\nu(x_0)$ of $x_0$. For every $\xi\in\Om$, we
denote by $[x_0,\xi[$ the geodesic ray in $T_\nu(x_0)$ starting from
$x_0$ and converging to $\xi$. We denote by $(x_n^\xi)_{n\in\NN}$ the
sequence of vertices of $[x_0,\xi[$, in this order along this ray. In
particular, $x_0^\xi=x_0$ and $d(x_k^\xi,x_n^\xi)=|k-n|$.

\smallskip\noindent \begin{minipage}{8.9cm} ~~~ 
Let $x\in VT_\nu(x_0)$. We define the {\it sector} of $x$ by
$$
\Om_x=\big\{\xi\in\Om\;:\; x\in[x_0,\xi[\big\}\;,
$$ 
the {\it cone} of $x$ by
$$
C_x=\big\{y\in VT_\nu(x_0)\;:\;
\exists\;\xi\in\Omega_x,\;\;y\in[x,\xi[\,\big\}\;,
$$
and, for every $n\in\NN$, the {\it sector-sphere} of $x$ of radius $n$
by
$$
S_x^n=C_x\cap S_{T_\nu(x_0)}(x_0,n)\;.
$$
\end{minipage}
\begin{minipage}{6cm}
\begin{center}
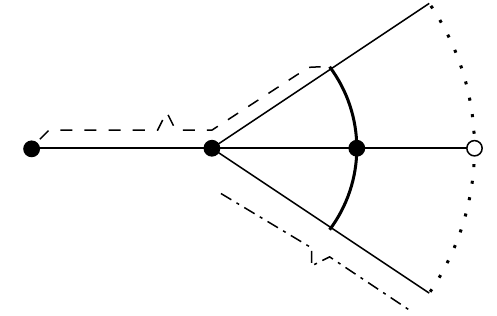
\end{center}
\end{minipage}

\medskip The {\it depth} of the cone $C_x$ or of the sector $\Om_x$ of
$x$ is defined to be the distance in the Hecke tree $T_\nu(x_0)$ from
$x$ to $x_0$. The sector-sphere $S_x^n$ is nonempty if and only if
$n$ is at least this depth. For every $\xi\in\Om$, the sequences
$(C_{x^\xi_n})_{n\in\NN}$ and $(\Om_{x^\xi_n})_{n\in\NN}$ are strictly
decreasing, with $\Om_{x_0}=\Om$, $C_{x_0}=VT_\nu(x_0)$,
$\bigcap_{n\in\NN} C_{x^\xi_n}=\emptyset$ and $\bigcap_{n\in\NN}
\Om_{x^\xi_n}=\{\xi\}$.

Note that if two cones (or sectors) intersect nontrivially, then one
of them is contained in the other. Also, sectors are nonempty
compact-open sets in $\Om$ and in particular contain infinitely many
rational ends, and even infinitely many $S$-rational ends.

\btheo\label{th:uncountable} There exists $c=c(x_0)>0$ such that, for every
$S$-cusp $z\in\E_\infty$ of $X_\infty$, the set of $\xi\in\Om$ having \cem\
towards the cusp $z$ is uncountable.  
\etheo

In particular, the set of $\xi\in\Om$ having \cem\ is uncountable.
Theorem \ref{u.cthm} in the Introduction follows immediately, being
the case when $\bf C$ is the projective line, in which case $X_\infty$
has only one cusp.

\medskip \dem Let $c=c(x_0)>0$ be the constant introduced in Theorem
\ref{theo:lossmass}. For every $S$-cusp $z\in\E_\infty$, we fix a
fundamental system $(V_n)_{n\in\NN}$ of open neighborhoods of $z$ in
$\wh{X_\infty}=X_\infty\cup\E_\infty$, so that $\{z\}=
\bigcap_{n\in\NN} V_n$. For all $n\in\NN$, let $\Sigma_n=\{0,1\}^n$ be
the set of words of length $n$ in $0$ and $1$. Let $\Sigma=
\bigcup_{n\in\NN} \Sigma_n$ be the set of finite words in $0$ and $1$.

We are going to define a map $\psi:\Sigma\to VT_\nu(x_0)$ with the
following properties: For all $n\in\NN$ and $\alpha\in\Sigma_n$,
\begin{enumerate}
\item[(1)] if $\beta$ is an initial subword of $\alpha$, then
  $\Om_{\psi(\alpha)}\subset \Om_{\psi(\beta)}$,
\item[(2)] if $\beta$ is an initial subword of $\alpha$ with
  $\beta\neq \alpha$, then the intersection $\Om_{\psi(\alpha 0)}\cap
  \Om_{\psi(\alpha 1)}$ is empty,
\item[(3)] the depth of the sector $\Om_{\psi(\alpha)}$ is at
  least $n$,
\item[(4)] we have $\mu_{\psi(\alpha)}(V_n)\geq c-\frac{1}{n+1}$.
\end{enumerate}

Assume for the moment that such a map $\psi$ is constructed. Let
$\Sigma_\infty=\{0,1\}^\NN$, which is uncountable. For every
$w\in\Sigma_\infty$, let $w_n$ be the initial subword of length $n$ of
$w$. Note that by Properties (1) and (3), for every
$w\in \Sigma_\infty$, the sequence of sectors
$(\Omega_{\psi(w_n)})_{n\in\NN}$ is strictly nested, and its
intersection contains a single point, denoted by $\xi_w$. Furthermore,
for every $n\in\NN$, we have $w_n\in [x_0, \xi_w[\,$. Note that by
Property (2), the map $w\mapsto \xi_w$ from $\Sigma_\infty$ to
$\Omega$ is injective. By Property (4), for every $w\in
\Sigma_\infty$, if $\theta_w$ is a weak-star accumulation point of
$(\mu_{\psi(w_n)})_{n\in\NN}$ in the space $\P(\widehat{X_\infty})$ of
probability measures on the compact space $X_\infty$, then
$\theta_w(\{z\})\geq c$.  Hence $\xi_w$ has \cem\ towards the cusp
$z$. This proves Theorem \ref{th:uncountable}.

\medskip We now build $\psi_{\mid \Sigma_n}$ by induction on
$n\in\NN$. Note that $\Sigma_0$ is reduced to the empty word
$\emptyset$, and define $\psi(\emptyset)=x_0$. Let $n\in\NN$, assume
that $\psi_{\mid \Sigma_n}$ is constructed, satisfying Properties
(1)--(4) for every $\alpha\in\Sigma_n$. For every $\alpha\in \Sigma_n$
and $j\in\{0,1\}$, let us define $\psi(\alpha j)$.

\begin{center}
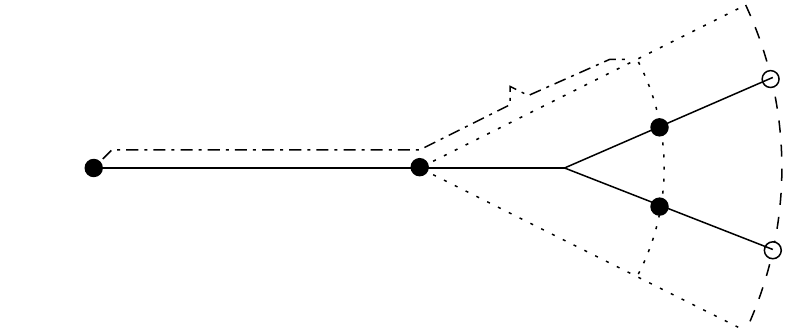
\end{center}

By density, there exist distinct points $\xi_0$ and $\xi_1$ in
$\Omega_{\psi(\alpha)}$ which are rational and whose associated cusps
$z_{\xi_0}$ and $z_{\xi_1}$ of $X_\infty$ respectively are both equal
to $z$.  By Theorem \ref{theo:lossmass}, $\xi_0$ and $\xi_1$ both have
\cem\ towards the cusp $z$. 

For all $j\in\{0,1\}$ and $m\geq d(x_0,\psi(\alpha)) +1$, the sector
$\Omega_{x_m^{\xi_j}}$ is strictly contained in
$\Omega_{\psi(\alpha)}$ and has depth at least $n+1$ by induction.
Since $\xi_0\neq \xi_1$, there exists $m_0\in\NN$ such that
$x_{m_0}^{\xi_0}\neq x_{m_0}^{\xi_1}$, so that for every $m,m'\geq
m_0$, the sectors $\Omega_{x_m^{\xi_0}}$ and $\Omega_{x_{m'}^{\xi_1}}$
are disjoints.

Let $j\in\{0,1\}$. We claim that there exists $n_j\geq m_0$ such that
$\mu_{x_{n_j}^{\xi_j}}(V_{n+1})\geq c-\frac{1}{n+2}$. Otherwise, for every
accumulation point $\theta$ of $\big(\mu_{x_m^{\xi_j}}\big)_{m\in\NN}$, we
have $\theta(\{z\})\leq c-\frac{1}{n+2}$, which contradicts the fact
that $\xi_j$ has \cem\ towards the cusp $z$.

Defining $\psi(\alpha 0)=x_{n_0}^{\xi_0}$ and $\psi(\alpha 1)=x_{n_1}^{\xi_1}$
gives the result.  
\cqfd

\subsection{Effective equidistribution of sector-spheres}
\label{subsec:effequisectorsphere} 

The aim of this section is to prove an effective statement regarding
the equidistribution in $X_\infty$ of the sector-spheres of the
vertices of the Hecke tree of $x_0$, Theorem \ref{theo:eed}, by using
the effective decay of matrix coefficients for the action of $G_S$ on
$\LL^2(X_S)$. This sectorial effective equidistribution result will be
the main tool used in Subsection \ref{subsec:exotic} in order to prove
Theorem \ref{nac} and its improvements. We first
introduce some notation.

\medskip We denote by $|E|$ the cardinality of any finite set $E$ and
by $\Delta_x$ the unit Dirac mass at any point $x$ of any measurable
space. For all $x\in VT_\nu(x_0)$ and $n\in\NN$ with $n\geq k$ where
$k=d_{T_\nu(x_0)}(x_0,x)$ is the depth of the sector $C_x$, let
$\eta_{n,\,x}$ be the uniform probability measure on the (finite
nonempty) sector-sphere $S_x^n$:
$$
\eta_{n,\,x}=\frac{1}{|S_x^n|}\;\sum_{y\in S_x^n}\;\Delta_y\;,
$$
that we consider as a probability measure on the locally compact space
$X_\infty$ with support $S_x^n$. Since the $\nu$-Hecke tree of $x_0$
(as is the Bruhat-Tits tree $\TT_\nu$) is $|\PP^1(k_\nu)|$-regular,
note that $|S_x^n|= |k_\nu|^{n-k}$ if $x\neq x_0$ and $n\geq k$, and that
$|S_x^n|= (|k_\nu|+1)|k_\nu|^{n-1}$ if $x= x_0$ and $n>0$.

For every place $\omega\in \P$, we define $W_\omega=\uG(O_\omega)$,
which is a maximal compact-open subgroup of $G_\omega$, and $W_S=
W_\infty\times W_\nu\subset G_\infty \times G_\nu=G_S$, which is a
maximal compact-open subgroup of $G_S$.

We denote by $m_\infty$ (respectively $m_S$) the Haar measure on
$G_\infty$ (respectively $G_S$), normalized so that
$m_\infty(W_\infty) =1$ (respectively $m_S(W_S)=1$). We again denote
by $m_\infty$ (respectively $m_S$) the measure on $X_\infty$
(respectively $X_S$) such that the covering map $G_\infty\ra X_\infty
=\Ga_\infty \bs G_\infty$ (respectively $G_S\ra X_S=\Ga_S\bs G_S$)
locally preserves the measures. Note that this measure on $X_\infty$
(respectively $X_S$) is nonzero and finite, but is not necessarily a
probability measure, the above normalisation of the Haar measures will
turn out to be more convenient. For every $k\in[1,+\infty]$, we define
$\LL^k(X_\infty)= \LL^k(X_\infty,m_\infty)$ (respectively
$\LL^k(X_S)=\LL^k(X_S,m_S)$~).

The group $G= G_\infty$ (respectively $G=G_S$) acts (on the left) on
the complex vector space of maps $\psi$ from $X=X_\infty$
(respectively $X=X_S$) to $\CC$, by right translation on the source:
For every $g\in G$, if $R_g:X\ra X$ is the right translation $x\mapsto
xg$, then $g\psi= \psi \circ R_g:x\mapsto \psi (xg)$.

A map $\psi$ from $X$ to $\CC$ is {\it locally constant} if there
exists a compact-open subgroup $U$ of $W=W_\infty$ (respectively
$W=W_S$) which leaves $\psi$ invariant:
$$
\forall\;g\in U,\;\;\; g\psi=\psi\;,
$$
or equivalently, if $\psi$ is constant on each orbit of $U$ under the
right action of $G$ on $X$. Note that $\psi$ is continuous, since the
orbits of $U$ are compact-open subsets.
We define
$$
d_f=\dim (\operatorname{Vect}_\CC W \psi)
$$
as the dimension of the complex vector space generated by the images
of $\psi$ under the elements of $W$, which is finite, and even
satisfies $d_f\leq [W:U]$. We define the {\it lc-norm} of every
bounded locally constant map $\psi:X\ra\CC$ by
$$
\|\psi\|_{lc}=\sqrt{d_\psi }\,\|\psi\|_\infty\;.
$$
Though the lc-norm does not satisfy the triangle inequality, we have
$\|\lambda\psi\|_{lc}=|\lambda|\,\|\psi\|_{lc}$ for every
$\lambda\in\CC$. We denote by $lc(X)$ the vector space of bounded
locally constant maps $\psi$ from $X$ to $\CC$.

Finally, given a set $A$ and maps $f,g:A\ra[0,+\infty[$ , we will
write $f\ll g$ if there exists a constant $c'>0$ such that $f(a)\leq
c' g(a)$ for all $a\in A$. If $f$ and $g$ depend on a parameter $p$,
we write $f\ll_p g$ if there exists a constant $c'>0$, possibly
depending on the parameter $p$, such that $f(a)\leq c' g(a)$ for all
$a\in A$.

\medskip The following result strenghtens the well-known result of
equidistribution of full Hecke spheres (see for instance the works of
Dani-Margulis \cite{DanMar93}, Clozel-Oh-Ullmo \cite{CloOhUll01},
Clozel-Ullmo \cite{CloUll04}, Eskin-Oh \cite{EskOh06b}, Benoist-Oh
\cite{BenOh07} in characteristic $0$), to an equidistribution result
of sector-spheres, which is furthermore effective. Taking $x=x_0$
gives as a particular case an effective equidistribution result of the
full Hecke spheres.

\btheo\label{theo:eed} There exists $\delta>0$ such that for every
$x\in VT_\nu(x_0)$, we have
\begin{equation}\label{eq:924}
\Big|\, \frac{m_\infty(\psi)}{m_\infty(X_\infty)}-
\eta_{n,\,x}(\psi) \,\Big|\ll \|\psi\|_{lc}\;  e^{-\delta\,n} 
\end{equation} 
for all $n\gg_x 1$ and $\psi\in lc(X_\infty)$.
\etheo

\dem Let us fix $x\in VT_\nu(x_0)$ and $\xi=\xi_x\in \Omega_x$, so
that $x=x^\xi_k$ for some fixed $k=k_x\in\NN$ (see the picture at the
beginning of Subsection \ref{subsec:uncountableeom}).

\medskip
\noindent{\bf Step 1: Thickening the sector-spheres. }
Note that the sector-spheres are measure zero subsets of
$X_\infty$. In order to be able to apply (effective) mixing arguments,
we have to replace them by (regular) bump fonctions around them. In
this step, we will define nice compact-open neighborhoods of the
sector-spheres, whose characteristic functions will be our bump
fonctions. By the construction of the sector-spheres, it is more
natural to lift the sector-spheres in $X_S$ and to work in the bundle
$X_S$ over $X_\infty$.

We will hence use a lot the $W_\nu$-bundle map $\pi_\infty$ (see
Subsection \ref{subsec:BruhatTitstrees}) from $X_S=\Ga_S\bs G_S$ to
$X_\infty=\Ga_\infty\bs G_\infty$, defined by $\Ga_S(g,h)\mapsto
\Ga_\infty g$ whenever $h\in W_\nu$. Recall (see Subsection
\ref{subsec:Hecketrees}) that the map $\hec_{g_0}$ from the
Bruhat-Tits tree $\TT_\nu$ to the Hecke tree $T_\nu(x_0)$, defined on
$V\TT_\nu= G_\nu/W_\nu$ by $h W_\nu\mapsto \pi_\infty(\Ga_S(g_0,h))$
is an isomorphism of trees, and we identify
$\partial_\infty\TT_\nu=\PP_1(K_\nu)$ and $\Omega$ by (the extension
to the boundary at infinity of) this map. We endow $T_\nu(x_0)$ with
the (left) action of $G_\nu$ making $\hec_{g_0}$ equivariant. Since
$W_\nu=\uG(O_\nu)$ acts transitively on $\Omega=\PP_1(O_\nu)$, we also
fix $w=w_x\in W_\nu$ such that $w\infty= \xi$, where $\infty=[1:0]$.

For all $n\in\NN$, we denote by $B^n_\nu$ the stabiliser in $W_\nu$ of
the point $x_n^\infty$ at distance $n$ from $x_0$ on the geodesic ray
$[x_0, +\infty[$ in the Hecke tree $\TT_\nu(x_0)$. The group
$B_\nu=B^k_\nu$ acts transitively on the sector-spheres
$S^n_{x_k^\infty}$ of $x_k^\infty$ for all $n\in\NN$.  As we have
already seen, for all $n\in\NN$, we have
$$
x_n^\infty=\hec_{g_0}(a_\nu^n*_\nu)=\pi_\infty(\Ga_S(g_0,a_\nu^n))\;.
$$
Note that $x^\xi_n=w x_n^\infty$ for all $n\in\NN$. In particular,
$x=w x_k^\infty$, hence $w B_\nu w^{-1}$ is the stabilizer in $W_\nu$
of $x$. It acts transitively on the sector-spheres $S_x^n$ of $x$ for
all $n\in\NN$, with stabilizer of $x^\xi_n$ equal to $w B^n_\nu
w^{-1}$.  Therefore, for all $n\in\NN$,
\begin{equation}\label{eq:1443}
S_x^n = w B_\nu w^{-1}x^\xi_n= w B_\nu x^\infty_n=
\pi_\infty(\Ga_S(g_0, w B_\nu a_\nu^n))\;.
\end{equation}
Now that we have this nice description of the sector-spheres, let us
define nice neighborhoods of them.

\blemm\label{lem:constructBeps} There exist $\sigma_1,\sigma_2>0$ and
a nondecreasing family $(B^\epsilon_\infty)_{\epsilon >0}$ of
compact-open subgroups of $W_\infty$, which is a fundamental system of
neighborhoods of the identity element in $W_\infty$, and which
satisfies
\begin{equation}\label{eq:1444}
\forall \;\epsilon>0,\;\;\;\;\;\; \sigma_1\;\epsilon^{-1}\leq
[W_\infty:B^\epsilon_\infty] \leq \epsilon^{-1}\;,
\end{equation}
and
\begin{equation}\label{eq:2444}
  \forall \;a\in A_\infty,\;\;\;\;\;\; 
a^{-1}B^\epsilon_\infty a \subset B^{\,\epsilon\, e^{\sigma_2|v_\infty(a)|}}_\infty\;.
\end{equation}
\elemm

\dem 
For every $n\in\NN$, let $Z_n$ be the kernel of the reduction
modulo $\pi_\infty^{n+1}$ map from $W_\infty=\uG(O_\infty)$ to the
finite group $\uG(O_\infty/ \pi_\infty^{n+1}O_\infty)$. Let us
consider $B^\epsilon_\infty = Z_{n_\epsilon}$ with $n_\epsilon =
\lfloor \frac{-\log(\epsilon[W:Z_0])} {\log\,[Z_1:Z_0]}\rfloor$ and
$\sigma_1= \frac{1}{[Z_1:Z_0]}$. Then $B^\epsilon_\infty$ is a
compact-open subgroup of $W_\infty$, we have $B^\epsilon_\infty\subset
B^{\epsilon'}_\infty$ if $\epsilon\leq\epsilon'$ and
$\bigcap_{\epsilon>0} B^\epsilon_\infty=\{1\}$. Equation
\eqref{eq:1444} follows since the index $[Z_{n+1}:Z_n]$ is constant,
hence $[W:Z_n]=[W:Z_0][Z_1:Z_0]^n$.

For all $a,b,c,d\in K_\infty$ and $t\in K_\infty^\times$, we have
$$
\begin{pmatrix} 1 & 0 \\ 0 & t^{-1}\end{pmatrix}
\begin{pmatrix} a & b \\ c & d \end{pmatrix}
\begin{pmatrix} 1 & 0 \\ 0 & t \end{pmatrix}
=\begin{pmatrix} a & t\,b \\ t^{-1}c & d \end{pmatrix}\;.
$$
Hence, using the isomorphism $\alpha_\infty:K_\infty^\times\ra
A_\infty$ defined just above Remark \ref{rem:extension}, we have
$a^{-1}Z_n a\subset Z_{n-|v_\infty(a)|}$ for all $a\in A_\infty$ and
$n\geq |v_\infty(a)|$ in $\NN$. Equation \eqref{eq:2444} (which will
only be used in Subsection \ref{subsec:exotic}) follows with
$\sigma_2=\log\, [Z_1:Z_0]$.  \cqfd

\medskip
For every $\epsilon>0$, we finally define the following compact-open
subset of $X_S$
$$
U_{\epsilon} = \Ga_S (g_0 B_\infty^\epsilon, w B_\nu)\;, 
$$
so that, for all $n\in\NN$, the image $\pi_\infty(U_{\epsilon}
a_\nu^n)$ of its translate by $a_\nu^n$ is a (small when $\epsilon$ is
small) neighborhood of the sector-sphere $S_x^n$ in $X_\infty$, by
Equation \eqref{eq:1443}.

\medskip
\noindent{\bf Step 2: Using the decay of matrix coefficients. }
In this step, we use the following theorem about effective decay of
matrix coefficients for the action of $G_S$ on $\LL^2(X_S)$ (see for
instance \cite{AthGhoPra12}). For every $g=(g_\infty,g_\nu)\in G_S=
G_\infty\times G_\nu$, we denote by $|g|_S$ the maximum of the norms
of the adjoint representations of $g_\infty$, $g_\nu$ (for the
operator norm on the $3\times 3$ matrices with entries in $K_\infty,
K_\nu$).

\btheo\label{theo:edmc} There exists $\delta_1>0$ such that 
\begin{equation}\label{eq:1102}
\big|\,m_S(\wt \psi\; \wt \varphi\circ R_g)- \frac{1}{m_S(X_S)}
m_S(\wt \psi)\,m_S(\wt \varphi)\,\big|
\ll \sqrt{d_{\wt \varphi}\, d_{\wt \psi}}\;
\|\wt \varphi\|_{2}\,\|\wt \psi\|_{2}\, |g|_S^{-\delta_1}
\end{equation}
for all locally constant maps $\wt \varphi,\wt \psi\in \LL^2(X_S)$ and
for every $g\in G_S$. \cqfd \etheo

Now, let us fix $\psi\in lc(X_\infty)$. We denote by $\wt \psi=
\psi\circ\pi_\infty$ its lift to $X_S$, which is constant on each
right $W_\nu$-orbit, hence is locally constant (since invariant under
$U\times W_\nu$ if $\psi$ is invariant under $U$). Note that $\wt
\psi\in \LL^2(X_S)$ since $m_S$ is finite and $\wt \psi$ is
bounded. By the normalization of the Haar measures, we have
$$
m_S(\wt \psi)=m_\infty(\psi)\;\;\;{\rm and}\;\;\;
m_S(X_S)=m_\infty(X_\infty)\;.
$$ 
Since $\sqrt{d_{\wt \psi}}=\sqrt{d_{\psi}}$ and $\|\wt \psi\|_{2}\leq
\sqrt{m_S(X_S)}\;\|\wt \psi\|_{\infty}$, we have
$$
\sqrt{d_{\wt \psi}}\;\|\wt \psi\|_{2}\ll \|\psi\|_{lc}\;.
$$ 
For every $\epsilon>0$, let $\varphi_\epsilon=
\frac{1}{m_S(U_\epsilon)}\mathbbm{1}_{U_\epsilon}$ be the normalized
characteristic function of $U_\epsilon$, so that
$m_S(\varphi_\epsilon) =1$.  The map $\varphi_\epsilon:X_S\ra\CC$ is
locally constant, since it is invariant under the right action of the
compact-open subgroup $B_\infty^\epsilon \times B_\nu$ of $W_S$.
%
%
We have
$$
d_{\varphi_\epsilon}=\dim\;\operatorname{Vect}_\CC W_S\,\varphi_\epsilon\leq 
[W_S:B_\infty^\epsilon \times B_\nu]=
[W_\infty:B_\infty^\epsilon] [W_\nu: B_\nu]\;.
$$
Since $W_\nu$ is compact and acts freely on each of its orbits on
$X_S$, there exists $\epsilon_0=\epsilon_0(x)>0$ such that if
$\epsilon\in\,]0,\epsilon_0]$, the map from $B_\infty^\epsilon \times
B_\nu$ to $X_S$ defined by $(g,h) \mapsto \Ga_S(g_0g,wh)$ is
injective, and measure preserving with image $U_\epsilon$. Hence, by
the normalization of the Haar measures, we have, for every $\epsilon
\in\,]0,\epsilon_0]$,
\begin{equation}\label{eq:1101} 
\|\varphi_\epsilon\|_2=m_S(U_\epsilon)^{-\frac{1}{2}}
=(m_\infty(B_\infty^\epsilon)\,m_\nu(B_\nu))^{-\frac{1}{2}}
=([W_\infty:B_\infty^\epsilon]\,[W_\nu:B_\nu])^{-\frac{1}{2}}\;.
\end{equation}
We therefore have $\sqrt{d_{\varphi_\epsilon}} \;
\|\varphi_\epsilon\|_2 \leq 1$.  Note that for every $n\in\NN$,
$$
|a_\nu^{-n}|_S=\max \{|a_\nu^{-n}|_\infty,|a_\nu^{-n}|_\nu\}
=\max \{|\pi_\nu^{\pm n}|_\infty,|\pi_\nu^{\pm n}|_\nu\}
=\max \{|k_\infty|^{n|v_\infty(\pi_\nu)|},|k_\nu|^{n}\}\;.
$$

Applying Equation \eqref{eq:1102} to the functions $\wt \psi$, $\wt
\varphi=\varphi_\epsilon$ and taking $g=a_\nu^{-n}$, we hence have,
with $\delta_2=\delta_1 \max\{|v_\infty(\pi_\nu)|\log|k_\infty|,
\log|k_\nu|\} >0$, for every $\epsilon\in\,]0,\epsilon_0]$,
\begin{equation}\label{eq:1103} 
\Big|\,\frac{1}{m_S(U_\epsilon)}\int_{U_\epsilon a_\nu^n} \wt \psi \;dm_S
- \frac{m_\infty(\psi)}{m_\infty(X_\infty)}\,\Big| 
    \ll \|\psi\|_{lc}\; e^{-\delta_2 \,n} \;.
\end{equation}

Let us now relate, for $\epsilon$ small enough, the above quantity
$\frac{1}{m_S(U_\epsilon)}\int_{U_\epsilon a_\nu^n} \wt\psi \;dm_S$ to
the average
$$
\eta_{n,\,x}(\psi)=\frac{1}{|S_x^n|}\;\sum_{y\in S_x^n}\;\psi(y)
$$ 
of $\psi$ on the sector-sphere $S^n_x$.

Let $w_1,\dots, w_\ell$ be representatives of the right cosets in
$B_\nu/B^n_\nu$, so that $B_\nu$ is a disjoint union $B_\nu=
\coprod_{i=1}^\ell w_i B^n_\nu$ and $m_\nu(B_\nu)= [B_\nu:B^n_\nu]\,
m_\nu(B^n_\nu)$.  By the transitivity properties seen in Step 1, the
map from $B_\nu/B^n_\nu$ to $S_x^n$ defined by $[h]\mapsto whw^{-1}
x_n^\xi$ is a bijection. For every $y\in S_x^n$, let $i_y\in\{1,
\dots, \ell\}$ be such that 
$$
y=ww_{i_y}w^{-1} x_n^\xi=w w_{i_y} x_n^\infty=
\hec_{g_0}(w w_{i_y} a_\nu^n\, *_\nu)=
\pi_\infty(g_0,w w_{i_y} a_\nu^n)\;.
$$
Let
$$
V_y=\Ga_S (g_0 B_\infty^\epsilon, w w_{i_y} B^n_\nu)
$$
so that $V_y\,a_\nu^n=\Ga_S (g_0 B_\infty^\epsilon, w w_{i_y}
a_\nu^n(a_\nu^{-n}B^n_\nu a_\nu^n))$. Note that $a_\nu^{-n}B^n_\nu
a_\nu^n$ is contained in $W_\nu$, since $B^n_\nu$ stabilizes
$x_n^\infty= a_\nu^n x_0$, hence the restriction of $\pi_\infty$ to
$V_y\,a_\nu^n$ has image $y B_\infty^\epsilon$ and its fibers are
orbits of $a_\nu^{-n}B^n_\nu a_\nu^n$.  For every $\epsilon \in\,
]0,\epsilon_0]$, since the map $(g,h) \mapsto \Ga_S(g_0g,wh)$ from
$B_\infty^\epsilon \times B_\nu$ to $X_S$ is injective, we hence have
$$
U_{\epsilon} = 
\coprod_{i=1}^\ell \Ga_S (g_0 B_\infty^\epsilon, w w_i B^n_\nu)
= \coprod_{y\in S_x^n} V_y\;.
$$
Therefore, for every $\epsilon\in\,]0,\epsilon_0]$, using Equation
\eqref{eq:1101} and by desintegration of $m_S$, we have
\begin{align}
\frac{1}{m_S(U_\epsilon)}\int_{U_\epsilon a_\nu^n} \wt\psi \;dm_S
& =\frac{1}{m_\infty(B_\infty^\epsilon)\,m_\nu(B_\nu)}
\sum_{y\in S_x^n} \int_{V_y a_\nu^n}\wt\psi \;dm_S
\nonumber\\
& =\frac{1}{m_\infty(B_\infty^\epsilon)\,m_\nu(B^n_\nu)\,|S_x^n|}
\sum_{y\in S_x^n} \int_{\pi_\infty(V_y a_\nu^n)}m_\nu(a_\nu^{-n}B^n_\nu a_\nu^n)\;
\psi \;dm_\infty 
\nonumber\\ & =
\frac{1}{|S_x^n|} \sum_{y\in S_x^n} \frac{1}{m_\infty(B_\infty^\epsilon)}
\int_{y B_\infty^\epsilon} \psi \;dm_\infty \;.
\label{eq:1104}
\end{align}

Define the {\it $\epsilon$-thin part} $X_\infty^\epsilon$ of
$X_\infty$ as the set of points $z\in X_\infty$ such that the map from
$B_\infty^\epsilon$ to $X_\infty$ defined by $h\mapsto zh$ is not
injective. Since $\psi$ is locally constant, there exists
$\epsilon_1=\epsilon_1(\psi)>0$ such that if
$\epsilon\in\;]0,\epsilon_1]$, then $\psi$ is
$B^\epsilon_\infty$-invariant. If $y\in S_x^n- (S_x^n\cap
X_\infty^\epsilon)$ and if $\epsilon\in\;]0,\epsilon_1]$, then
\begin{equation}\label{eq:11045} 
\frac{1}{m_\infty(B_\infty^\epsilon)} \int_{y B_\infty^\epsilon} \psi
\;dm_\infty =\psi(y)\;.
\end{equation}
A trivial majoration gives
\begin{equation}\label{eq:1105} 
\Big|\,\frac{1}{|S_x^n|}\sum_{y\in S_x^n\cap X_\infty^\epsilon}
\Big( \psi(y)- \frac{1}{m_\infty(B_\infty^\epsilon)} \int_{B_\infty^\epsilon} 
\psi dm_\infty\Big)\,\Big|\leq 2\,\|\psi\|_\infty 
\frac{|S_x^n\cap X_\infty^\epsilon|}{|S_x^n|}\;.
\end{equation}

Separating the summation over $S_x^n$ on one hand over $S_x^n\cap
X_\infty^\epsilon$ and on the other hand over $S_x^n- (S_x^n\cap
X_\infty^\epsilon)$, for every $\epsilon\in\; ]0,\min\{\epsilon_0,
\epsilon_1\}]$, by Equations \eqref{eq:1103}, \eqref{eq:1104},
\eqref{eq:11045} and \eqref{eq:1105}, we hence have
\begin{equation}\label{eq:1106}
\Big|\,\eta_{n,\,x}(\psi)- \frac{m_\infty(\psi)}{m_\infty(X_\infty)}\,\Big| 
\ll \|\psi\|_{lc}\; e^{-\delta_2 \,n}+ \|\psi\|_\infty 
\frac{|S_x^n\cap X_\infty^\epsilon|}{|S_x^n|}\;.
\end{equation}


\medskip
\noindent{\bf Step 3: Estimating the thin part of sector-spheres. }
The aim of this step is to prove that the part of the sector-spheres
contained in the thin part of $X_\infty$ is negligible, if $\epsilon$
is well-chosen. More precisely, let us prove that there exists
$\delta_4>0$ such that for every $n\gg_x1$, if $\epsilon=
e^{-\delta_2n}$, then
\begin{equation}\label{eq:1110} 
\frac{|S_x^n \cap X_\infty^\epsilon|}{|S_x^n|} \ll  e^{-\delta_4n}
\end{equation}
for every $n\in\NN$ with $n>k$.

For this, we will apply the arguments of Step 2 to a particular map
$\psi=\psi_\epsilon$, where $\psi_\epsilon$ is, for every
$\epsilon>0$, the characteristic function of the $\epsilon$-thick part
$X_\infty-X_\infty^\epsilon$ of $X_\infty$. Note that $\psi_\epsilon$
is invariant under $B_\infty^\epsilon$, hence $\psi_\epsilon$ is
bounded and locally constant, and $\epsilon_1(\psi_\epsilon)=+\infty$.
Denoting by $\wt \psi_\epsilon$ the lift of $\psi_\epsilon$ to $X_S$,
by Equations \eqref{eq:1104} and \eqref{eq:11045} applied with
$\psi=\psi_\epsilon$, for every $\epsilon\in\;]0,\epsilon_0]$, we have
\begin{equation}\label{eq:1111} 
\frac{1}{m_S(U_\epsilon)}\int_{U_\epsilon a_\nu^n} \wt\psi_\epsilon \;dm_S 
=\frac{1}{|S_x^n|} \sum_{y\in S_x^n} \psi_\epsilon(y)=
\frac{|S_x^n-(S_x^n\cap X_\infty^\epsilon)|}{|S_x^n|} \;.
\end{equation}

By Equation \eqref{eq:1444}, we have
\begin{equation}\label{eq:1112} 
\|\psi_\epsilon\|_{lc}=\sqrt{d_{\psi_\epsilon}}\;\|\psi_\epsilon\|_{\infty}
\leq \sqrt{[W_\infty:B_\infty^\epsilon]}\leq \epsilon^{-\frac{1}{2}}\;.
\end{equation}

By the exponential decay of the volumes in the cusps of the graph of
groups $\Ga_\infty\bs\!\bs\TT_\infty$, hence in $X_\infty$, there exists
$\delta_3>0$ such that 
\begin{equation}\label{eq:1113}
m_\infty(X_\infty^\epsilon)\ll \epsilon^{\delta_3}
\end{equation}
for every $\epsilon>0$. 

For every $n\in\NN$, define $\epsilon=e^{-\delta_2n}$. Note
that if $n\gg_x 1$, then $\epsilon\leq \epsilon_0$ (recall that
$\epsilon_0$ depends on $x$). Therefore, using 

$\bullet$~ Equation \eqref{eq:1113} for the second inequality,

$\bullet$~ Equation \eqref{eq:1111} and the definition of
$\psi_\epsilon$ for the third line,

$\bullet$~ Equation \eqref{eq:1103} with $\psi=\psi_\epsilon$ for the
fourth inequality,

$\bullet$~ Equation \eqref{eq:1112} for the
fifth inequality,

$\bullet$~ the definition of $\epsilon$ and the constant
$\delta_4=\delta_2\min\{\delta_3,\frac{1}{2}\}>0$ for the last
inequality,

\noindent
we have, for all $n\gg_x 1$,
\begin{align*}
\frac{|S_x^n \cap X_\infty^\epsilon|}{|S_x^n|} &
\leq  
\Big|\frac{|S_x^n \cap X_\infty^\epsilon|}{|S_x^n|}-
\frac{m_\infty(X_\infty^\epsilon)}{m_\infty(X_\infty)}\Big|+
\frac{m_\infty(X_\infty^\epsilon)}{m_\infty(X_\infty)}\\
&\ll \Big|\frac{|S_x^n -(S_x^n \cap X_\infty^\epsilon)|}{|S_x^n|}-
\frac{m_\infty(X_\infty-X_\infty^\epsilon)}{m_\infty(X_\infty)}\Big|
+\epsilon^{\delta_3}\\
&=\Big|\frac{1}{m_S(U_\epsilon)}
\int_{U_\epsilon a_\nu^n} \wt\psi_\epsilon \;dm_S -
\frac{m_\infty(\psi_\epsilon)}{m_\infty(X_\infty)}\Big|+
\epsilon^{\delta_3}\\
&\ll \|\psi_\epsilon\|_{lc}\; e^{-\delta_2 \,n}+\epsilon^{\delta_3} 
\leq \epsilon^{-\frac{1}{2}}e^{-\delta_2 \,n}+\epsilon^{\delta_3} 
\leq  2\;e^{-\delta_4n}\;.
\end{align*}
This proves the claim \eqref{eq:1110} of Step 3.

\medskip
\noindent{\bf Step 4: Conclusion. }
Since $\|\psi\|_\infty \leq \|\psi\|_{lc}$, Theorem \ref{theo:eed} now
follows from Equations \eqref{eq:1106} and \eqref{eq:1110}, with
$\delta=\min\{\delta_2,\delta_4\}$.  
\cqfd

\subsection{Exotic behavior of $A_\infty$-periodic 
measures along Hecke rays}
\label{subsec:exotic} 

In this final Subsection, we use the tools introduced in Subsections
\ref{subsec:ratheckeray} and \ref{subsec:effequisectorsphere} to
construct even more exotic asymptotic behaviors of the
$A_\infty$-periodic measures $\mu_x$ as $x$ varies along geodesic rays
in the Hecke tree of $x_0$.

\btheo\label{th:general} Let $(\mu_i)_{i\in\NN}$ be an enumeration of
all periodic $A_\infty$-invariant probability measures on $X_\infty$,
and let $z\in\E_\infty$ be a $S$-cusp of $X_\infty$. There exist
$c,c'>0$ such that the set of $\xi\in\Om$ having $c$-escape of mass
towards $z$ and verifying
$$
\forall\; i\in\NN,\; \exists\;\theta_i\in\Theta_\xi,\;\;\;\;
c'\mu_i\leq\theta_i
$$
is uncountable. In particular, $\{\xi\in\Omega\;:\; |\Theta_\xi|=
\infty\}$ is uncountable.  
\etheo

Note that there are indeed only countably many periodic
$A_\infty$-orbits, and that this result immediately implies Theorem
\ref{nac}.

\medskip The proof of Theorem \ref{th:general} relies on the following
two lemmas. We consider again the family
$(B^\epsilon_\infty)_{\epsilon >0}$ of compact-open subgroups of
$W_\infty$ constructed in Lemma \ref{lem:constructBeps}.

\blemm\label{lem:close} There exists $\delta'>0$ such that for every
$x\in V T_\nu(x_0)$, for every $A_\infty$-periodic point $y_0\in
X_\infty$, and for every $n\gg_{x,y_0} 1$, the intersection $S_x^n
\cap (y_0 A_\infty B_\infty^{e^{-\delta'n}})$ is nonempty.  
\elemm

\dem Let $\delta'=\frac{\delta}{3}$ where $\delta$ is the constant
given by our effective equidistribution result of sector-spheres,
Theorem \ref{theo:eed}. For every $\epsilon>0$, let
$\psi_\epsilon=\mathbbm 1_{y_0 A_\infty B_\infty^\epsilon}$ be the
characteristic function of the $B_\infty^\epsilon$-thickening of the
periodic orbit $y_0 A_\infty$, which is bounded and locally constant.
We are going to use Theorem \ref{theo:eed} applied to
$\psi=\psi_\epsilon$ for a suitably choosen $\epsilon$.

There exists $\epsilon_2=\epsilon_2(y_0)>0$ such that if
$\epsilon\in\;]0,\epsilon_2]$, then the orbit map
$B_\infty^\epsilon\ra y_0 B_\infty^\epsilon$ is injective. Hence by
the normalisation of the Haar measure and by Equation \eqref{eq:1444},
we have, for every $\epsilon\in\;]0,\epsilon_2]$,
$$
m_\infty(\psi_\epsilon)=m_\infty(y_0 A_\infty B_\infty^\epsilon)
\geq m_\infty(y_0 B_\infty^\epsilon)=m_\infty(B_\infty^\epsilon)
=\frac{1}{[W:B_\infty^\epsilon]}\geq \epsilon\;.
$$
Furthermore
$$
\|\psi_\epsilon\|_{lc}=\sqrt{d_{\psi_\epsilon}}\;\|\psi_\epsilon\|_{\infty}
\leq \sqrt{[W:B_\infty^\epsilon]}\leq \epsilon^{-\frac{1}{2}}\;.
$$
By Theorem \ref{theo:eed}, there exists $\kappa>0$ such that if $n\gg_x 1$,
we have, for every $\epsilon\in\;]0,\epsilon_2]$,
$$
\eta_{n,\,x}(\psi_\epsilon)\geq 
\frac{m_\infty(\psi_\epsilon)}{m_\infty(X_\infty)}-
\kappa\;\|\psi_\epsilon\|_{lc}\; e^{-\delta\,n} 
\geq \frac{\epsilon}{m_\infty(X_\infty)}
-\kappa\;\epsilon^{-\frac{1}{2}}\;e^{-\delta\,n}\;.
$$
Let us now consider $\epsilon= 2 (\kappa\,m_\infty(X_\infty)\,
e^{-\delta\,n})^{\frac{2}{3}}$.  By the definition of $\delta'$, we have
$\epsilon\leq e^{-\delta'n}$ if $n\gg 1$. If $n\gg_{y_0} 1$, then
$\epsilon$ belongs to $]0,\epsilon_2]$. The previous centered formula
then gives $\eta_{n,\,x}(\psi_\epsilon) >0$ if $n\gg_{x,y_0} 1$. Hence
if $n\gg_{x,y_0} 1$, the support of the measure $\eta_{n,\,x}$, which is
$S^n_x$, meets the support of the function $\psi_\epsilon$, which is
$y_0 A_\infty B_\infty^\epsilon\subset y_0 A_\infty
B_\infty^{e^{-\delta'n}}$, as wanted.  \cqfd

\blemm\label{lem:macc} Let $y$ and $x_k$, for $k\in\NN$, be
$A_\infty$-periodic points in $X_\infty$. Suppose that there exist
$\sigma,\delta'>0$ and an increasing sequence $(n_k)_{k\in\NN}$ of
positive integers such that:
\begin{enumerate}
\item[(1)] For every $k\in\NN$, the period, under the geodesic
  flow in $\Ga_\infty\bs\G\TT_\infty$, of $\pi'_\infty(x_k)$ is at
  most $\sigma \,n_k$.
\item[(2)] There are infinitely many $k\in\NN$ such that
  $x_k \in y A_\infty B_\infty^{e^{-\delta'n_k}}$.
\end{enumerate}
Then there exists a weak-star accumulation point $\theta$ of
$(\mu_{x_k})_{k\in\NN}$ such that $\frac{\delta'}{2\,\sigma_2\,\sigma}\;\mu_y \leq
\theta$. 
\elemm

\dem Up to extracting a subsequence, we may assume that $x_k \in y
A_\infty B_\infty^{e^{-\delta'n_k}}$ for every $k\in\NN$.

By Assumption (1) and by the equivariance of the canonical bundle map
$\pi'_\infty:X_\infty=\Ga_\infty\bs G_\infty \ra \Ga_\infty\bs\G
\TT_\infty= \Ga_\infty\bs G_\infty/\uA(O_\infty)$ with respect to the
epimorphism $v_\infty: A_\infty\ra\ZZ$ where $\ZZ$ acts by the
geodesic flow (see Subsection \ref{subsec:BruhatTitstrees} and in
particular Equation \eqref{eq:equivarvomega}), we have
$$
x_k A_\infty= \{x_ka\;:\; a\in A_\infty\;\;{\rm and}\;\;
|v_\infty(a)| \leq \sigma \,n_k\}\;.
$$
For every $a\in A_\infty$ such that $|v_\infty(a)| \leq
\frac{\delta'}{2\sigma_2} \, n_k$, we have, by Assumption (2) and by
Equation \eqref{eq:2444},
$$
x_ka\in y A_\infty a^{-1}B_\infty^{e^{-\delta'n_k}}a\subset 
y A_\infty B_\infty^{e^{-\delta'n_k}e^{\sigma_2|v_\infty(a)|}}
\subset 
y A_\infty B_\infty^{e^{-\frac{\delta'}{2}\,n_k}}\;.
$$
If $m_{A_\infty}$ is the Haar measure of $A_\infty$ normalized so that
$\mu_{A_\infty}(\uA(O_\infty))=1$, then the pushforward of
$m_{A_\infty}$ by $v_\infty$ is the counting measure of $\ZZ$ and
$\mu_{x_n}$ is the measure on $x_kA_\infty$ induced by $m_{A_\infty}$,
normalized to be a probability measure. Hence $\{x_ka\;:\;
|v_\infty(a)| \leq \frac{\delta'}{2\,\sigma_2} \,n_k\}$ occupies at
least $\frac{\delta'}{2\,\sigma_2\,\sigma}$ of the total mass of $x_k
A_\infty$, and accumulates on $y A_\infty$. Hence at least
$\frac{\delta'}{2\,\sigma_2\,\sigma}$ of the total mass of any weak-star
accumulation point of $(\mu_{x_k})_{k\in\NN}$ is accumulated on
$\mu_y$. \cqfd

\medskip
\noindent
{\bf Proof of Theorem \ref{th:general}. } 
Let $(\mu_i)_{i\in\NN}$ and $z$ be as in this statement.

Let us denote by $(\eta_n)_{n\in\NN}$ a sequence of measures on
$X_\infty$ which contains the zero measure as well as all the
$A_\infty$-invariant probability measures of the $A_\infty$-periodic
points of $X_\infty$, in such a way that each measure appears
infinitely many times.  Using Lemma~\ref{lem:close} and arguing
similarly to the proof of Theorem~\ref{th:uncountable}, with
$(V_n)_{n\in\NN}$ a fundamental system of open neighborhoods of $z$ in
$\wh{X_\infty}=X_\infty\cup\E_\infty$, we can build inductively
uncountably many sequences $(x_k)_{k\in\NN}$ in $VT_\nu(x_0)$ such
that the following holds.
\begin{enumerate}
\item[(1)] The sequence of cones $(C_{x_k})_{k\in\NN}$ is strictly nested,
  so that if $\bigcap_{k\in\NN}\Omega_{x_k}=\{\xi\}$, then
  $(x_k)_{k\in\NN}$ is a subsequence of the sequence
  $(x_n^\xi)_{n\in\NN}$ of vertices of the Hecke ray from $x_0$ to $\xi$.
\item[(2)] If $(x_k)_{k\in\NN}\neq(x'_k)_{k\in\NN}$ are two of these
  sequences, then the sectors $\Omega_{x_k}$ and $\Omega_{x'_k}$ are
  disjoint for $k$ big enough. In particular, the map $(x_k)_{k\in\NN}
  \mapsto \xi=\displaystyle\lim_{k\ra\infty} x_k$ is injective, and
  there are uncountably many such $\xi$'s.
\item[(3)] For every $k\in\NN$, denoting by $n_k$ the depth of $x_k$ which
  we may assume to be at least $1$,
\begin{enumerate}
\item if $\eta_k=0$ then $\mu_{x_k}(V_k)\geq c-\frac{1}{k+1}$, where
  $c=c(x_0)>0$ is the constant introduced in Theorem
  \ref{theo:lossmass},
\item if $\eta_k$ is the $A_\infty$-invariant probability measure on
  the orbit of an $A_\infty$-periodic point $y_k$, then $x_k\in y_k
  A_\infty B_\infty^{e^{-\delta'n_k}}$.
\end{enumerate}
\end{enumerate}
Since Case (a) occurs infinitely many times, the set $\Theta_\xi$
contains a weak-star accumulation point $\theta$ of
$(\mu_{x_k})_{k\in\NN}$ such that $\theta(\{z\})\geq c$.

Let $i\in\NN$, and let $y_i$ be in the support of $\mu_i$. By Case
(b), since there are infinitely many $k\in\NN$ such that $\eta_k=
\mu_i$, there are infinitely many $k\in\NN$ such that $x_k\in y_i
A_\infty B_\infty^{e^{-\delta'n_k}}$. With the terminology of
Subsection \ref{subsec:periodicorbit}, we use Theorem
\ref{theo:lineargrowpspherorbit} applied to a loxodromic element
$\ga_0$ associated to the choosen representative $g_0$ of the
$A_\infty$-periodic point $x_0$. This result gives that the period,
under the geodesic flow in $\Ga_\infty\bs\G\TT_\infty$, of
$\pi'_\infty(x_k)$ (since $x_k$ has depth $n_k\geq 1$ in $T_\nu(x_0)$)
is at most $\sigma \,n_k$ for some $\sigma>0$ (depending only on
$p,\ga_0,\nu$). Applying Lemma~\ref{lem:macc} with $y=y_i$, the set
$\Theta_\xi$ contains a weak-star accumulation point $\theta_i$ of
$(\mu_{x_k})_{k\in\NN}$ such that $\frac{\delta'}{2\,\sigma\,\sigma_2}
\;\mu_i\leq \theta_i$.

This proves the result, with $c'=\frac{\delta'}{2\,\sigma\,\sigma_2}$
(which does not depend on $i$). \cqfd

 {\small \bibliography{../biblio} }

\bigskip
{\small
\noindent \begin{tabular}{l} 
Mathematics Department, Technion\\
Israel Institute of Technology, Haifa, 32000 ISRAEL.\\
{\it e-mail:  alikkimron@gmail.com}
\end{tabular}
\medskip

\noindent \begin{tabular}{l}
Laboratoire de math\'ematique d'Orsay,
\\ UMR 8628 Univ. Paris-Sud, CNRS\\
Universit\'e Paris-Saclay,
91405 ORSAY Cedex, FRANCE\\
{\it e-mail: frederic.paulin@math.u-psud.fr}
\end{tabular}
\medskip

\noindent \begin{tabular}{l} 
Mathematics Department, Technion \\
Israel Institute of Technology, Haifa, 32000 ISRAEL.\\
{\it e-mail: ushapira@tx.technion.ac.il}
\end{tabular}
}

\end{document}